\documentclass[reqno]{amsart}
\usepackage{hyperref}
\usepackage{amssymb,amsmath,amsthm, enumerate}
\usepackage{xcolor}
\usepackage[pdftex]{graphics}
\newcommand{\R}{\mathbb{R}}
\newcommand{\Z}{\mathbb{Z}}
\newcommand{\N}{\mathbb{N}}

\def\R{{\mathbb{R}}}

\def\N{{\mathbb{N}}}
\def\Z{{\mathbb{Z}}}
\def\F{{\mathbb{F}}}

\providecommand{\keywords}[1]
{
  \small	
  \textbf{\textit{Keywords---}} #1
}

\begin{document}


\title{Multiplicity results for non-local operators of elliptic type}



\author{Emer Lopera}
\address{Universidad Nacional de Colombia-Manizales, Colombia.}
\curraddr{}
\email{edloperar@unal.edu.co}

\author{Leandro Rec\^{o}va}
\address{California State Polytechnic University, Pomona - USA.}
\curraddr{}
\email{llrecova@cpp.edu}

\author{Adolfo Rumbos}
\address{ Pomona College, Claremont, California - USA.}
\curraddr{610 N. College Avenue, 91711 Claremont, CA}
\email{arumbos@pomona.edu}
\keywords{Mountain Pass Theorem, Morse Theory, Critical Groups, local linking.}

\subjclass[2010]{Primary 35J20 }


\begin{abstract}
  In this paper, we study a class of problems proposed by Servadei and Valdinoci in \cite{Ser3}; namely,
   \begin{equation}\label{prob_0}
   \left\{\begin{aligned} -\mathcal{L}_{K} u(x)-\lambda u(x) & =f(x,u), \mbox{ for } x\in \Omega; \\ 
    u & =0 \quad \text{ in } \R^{N}\backslash\Omega,
    \end{aligned} \right.
\end{equation}
where $\Omega\subset \R^{N}$ is an open bounded set with Lipschitz boundary, $\lambda\in\R$, $f\in C^{1}(\overline{\Omega}\times\R,\R)$, with $f(x,0) = 0$ 
for $x\in\Omega$, and $\mathcal{L}_K$ is a non-local integrodifferential operator with homogeneous Dirichlet boundary condition. By computing the critical groups of the associated energy functional for problem (\ref{prob_0}) at the origin and at infinity, respectively, we prove that problem (\ref{prob_0}) has three nontrivial solutions for the case $\lambda < \lambda_1$ and two nontrivial solutions for the case $\lambda\geqslant\lambda_1,$ where $\lambda_1$ is the first eigenvalue of the operator $-\mathcal{L}_K$. Finally, assuming that the nonlinearity $f$ is odd in the second variable, we prove the existence of an unbounded sequence of weak solutions of problem (\ref{prob_0}) for the case
 $\lambda\geqslant\lambda_1$. We use variational methods and infinite-dimensional Morse theory to obtain the results. 
  
\end{abstract}
\numberwithin{equation}{section}
\newtheorem{theorem}{Theorem}[section]
\newtheorem{lemma}[theorem]{Lemma}
\newtheorem{definition}[theorem]{Definition}
\newtheorem{proposition}[theorem]{Proposition}
\newtheorem{prop}[theorem]{Proposition}
\newtheorem{corollary}[theorem]{Corollary}
\newtheorem{remark}[theorem]{Remark}
\allowdisplaybreaks
\maketitle

\section{Introduction}\label{secint}
 
In this paper, we study a class of problems involving non-local, integro-differential operators with homogeneous Dirichlet boundary conditions proposed by Servadei and Valdinoci in \cite{Ser3} and Servadei in \cite{Ser4}.  We start by introducing some notation before presenting the main results.

Let $K:\R^{N}\backslash\{0\}\rightarrow (0,\infty)$ be a function that satisfies the following assumptions:
 \begin{itemize}
     \item [($H_1$)] $\gamma K\in L^{1}(\R^{N})$ where $\gamma(x)=\min\{|x|^2,1\}$, for $x\in\R^N$. 
     \item [($H_2$)] There exist $s\in(0,1)$  and 
     $\theta > 0$  such that 
     $$K(x)\geqslant \theta|x|^{-(N+2s)},\quad\mbox{ for }x\in\R^{N}\backslash\{0\}.$$
     \item [$(H_3)$] $K(x)=K(-x)$, for any $x\in\R^{N}\backslash\{0\}$. 
 \end{itemize}

We recall that  a set $\Omega\subset\R^{N}$ has a Lipschitz boundary, and call it a Lipschitz domain, if, for every $x_0\in\partial\Omega$, there exists $r>0$ and a map $h:B_{r}(x_0)\rightarrow B_{1}(0)$ such that
 \begin{enumerate}
     \item [(i)] $h$ is a bijection;
     \item [(ii)] $h$ and $h^{-1}$ are both Lipschitz continuous functions;
     \item [(iii)] $h(\partial\Omega\cap B_{r}(x_0))=Q_0$;
     \item [(iv)] $h(\Omega\cap B_{r}(x_0))=Q_{+}$,
 \end{enumerate}
where $B_{r}(x_0)$ denotes the $N$--dimensional open ball of radius $r$ and center at $x_0\in\partial\Omega$, and 
$$Q_{0}:=\{(x_1,\ldots,x_N)\in B_{1}(0)\mbox{ }| x_N=0\}
$$
and
$$
Q_+:=\{(x_1,\ldots,x_N)\in B_1(0)\mbox{ }|x_N > 0\}.
$$

Denote by $\Omega\subset\R^{N}$ an open bounded set with Lipschitz boundary.  Fix $s\in (0,1)$ such that
$N>2s$, and define $X$ to be the linear space of Lebesgue measurable functions from $\R^{N}$ to $\R$ such that their restrictions to $\Omega$ belong to $L^2(\Omega)$ and the map
$(x,y)\mapsto ((u(x)-u(y))\sqrt{K(x-y)}$ is in
$L^{2}(\R^{2N}\backslash(\Omega^{c}\times\Omega^{c}))$
for all $u\in X$.
We endow $X$ with the norm
$\|\cdot\|_X\colon X\to\mathbb{R}$ given by 
\begin{equation}
\|u\|^2_X := \|u\|^2_{L^2(\Omega)} + \int_Q |u(x)-u(y)|^2K(x-y)\,dx\,dy ,
\end{equation}
for $u\in X$,
where $Q=\R^{2N}\backslash (\Omega^c\times\Omega^c).$

Define $X_0$ to be
\begin{equation}X_0=\{u\in X\;|\;u=0\mbox{ a.e in }\R^{N}\backslash\Omega\}.
\label{x0def}
\end{equation}
It was proved in \cite[Lemma 7]{Ser1} that the space $X_0$ defined in  (\ref{x0def}) is a Hilbert space 
with inner product given by 
\begin{equation}
    \langle u,v\rangle _{X_0} = \int_{Q}(u(x)-u(y))(v(x)-v(y))K(x-y)\,dx\,dy,
    \label{scalarproddef}
\end{equation}
for $u,v\in X_{0}.$ 
Thus, a norm in $X_0$ is given by  
\begin{equation}
    \|u\|_{X_0} = \sqrt{\langle u,u\rangle_{X_0}}=\left(\int_{Q}|u(x)-u(y)|^2K(x-y)\,dx\,dy\right)^{1/2}, 
    \label{normx0}
\end{equation}
for $u\in X_0$.

By virtue of condition $(H_1)$, it can be shown that $C_{0}^{2}(\Omega)\subset X_{0}$ which implies that $X_0$ and $X$ are non-empty. For a proof of this result, see \cite[Lemma $7$]{Ser5}.  

The spaces $X$ and $X_0$ are related to fractional Sobolev spaces as we will briefly describe below. 

For $s\in (0,1)$ and $p\in [1,+\infty)$, define the fractional Sobolev space $W^{s,p}(\Omega)$ as follows:
$$W^{s,p}(\Omega)=\left\{u\in L^{p}(\Omega):\frac{\left|u(x)-u(y)\right|}{\left|x-y\right|^{\frac{N}{p}+s}}\in L^{p}(\Omega\times\Omega)\right\},$$
endowed with the norm 
$$\|u\|_{W^{s,p}(\Omega)}=\left(\int_{\Omega}|u|^{p}\,dx + \iint_{\Omega}\frac{\left|u(x)-u(y)\right|^{p}}{\left|x-y\right|^{N+sp}}\,dx\,dy \right)^{\frac{1}{p}},
    \mbox{ for }
        u\in W^{s,p}(\Omega).
$$

In the case $p=2$, we write $H^{s}(\Omega)=W^{s,2}(\Omega)$, and note that $H^s(\Omega)$ is a Hilbert 
space.  The following lemma provides the connection between the spaces $X$ and $X_{0}$ and $H^{s}(\Omega)$. 

\begin{lemma}\label{SVLem5}\cite[Lemma $5$]{Ser1}
Let $K:\R^{N}\backslash\{0\}\rightarrow (0,+\infty)$ satisfy assumptions $(H_1)$--$(H_3)$. Then, the following assertions hold true:
\begin{itemize}
    \item [(a)] if $u\in X$, then $u\in H^{s}(\Omega)$. Moreover
    $$\|u\|_{H^{s}(\Omega)}\leqslant c(\theta)\|u\|_{X};$$
    \item [(b)] if $u\in X_{0}$, then $u\in H^{s}(\R^{N}).$ Moreover
    $$\|u\|_{H^{s}(\Omega)}\leqslant \|u\|_{H^{s}(\R^{N})}\leqslant c(\theta)\|u\|_{X}.$$
\end{itemize}
    In both cases $c(\theta)=\max\{1,\theta^{-1/2}\},$ where $\theta$ is given in $(H_2)$. 
\end{lemma}

We will also need a result about fractional Sobolev embeddings in this work. We start with the definition of an extension domain.

\begin{definition}\cite[Section $5$, page $545$]{Nezza1}
    For any $s\in(0,1)$ and any $p\in[1,+\infty)$, we say that an open set $\Omega\subset\R^{N}$ is an extension domain for $W^{s,p}$ if there exists a positive constant $C=C(N,p,s,\Omega)$ such that: for every function $u\in W^{s,p}(\Omega)$ there exists $\widetilde{u}\in W^{s,p}(\R^{N})$ with $\widetilde{u}=u$ for all $x\in\Omega$ and 
    $$\|\widetilde{u}\|_{W^{s,p}(\R^{N})}\leqslant C\|u\|_{W^{s,p}(\Omega)}.$$
\end{definition}
We note that if $\Omega\subset\R^{N}$ is bounded with Lipschitz boundary, then $\Omega$ is an extension domain for $W^{s,p}$ (see Theorem 5.2 in
\cite{Nezza1}).

The next theorem provides the Sobolev inequality for fractional Sobolev spaces. 

\begin{theorem}\cite[Theorem $6.7$]{Nezza1}
Let $s\in(0,1)$ and $p\in [1,+\infty)$ be such that $sp < N$. Let $\Omega\subset\R^{N}$ be an extension domain for $W^{s,p}$. Then there exists a positive constant $C=C(N,p,s,\Omega)$ such that, for any $u\in W^{s,p}(\Omega)$, we have 
\begin{equation}
    \|u\|_{L^{r}(\Omega)}\leqslant C\|u\|_{W^{s,p}(\Omega)},
    \label{sobineq}
\end{equation}
for any $r\in [p,p^{*}]$, where $p^{*} = \dfrac{Np}{N-sp}$ is the fractional critical exponent; 
i.e, the space $W^{s,p}(\Omega)$ is continuously embedded in $L^{r}(\Omega)$ for any $r\in[p,p^{*}].$ If, in addition, $\Omega$ is bounded, then the space $W^{s,p}(\Omega)$ is continuously embedded in $L^{r}(\Omega)$ for any $r\in [1,p^{*}]$. 
\label{theosobembed}
\end{theorem}

As a consequence of the results in Lemma \ref{SVLem5} and the embedding in Theorem
\ref{theosobembed}, we have the following embedding
result for $X_0$:
\begin{proposition}\label{X0EmbedThm}
Let $K\colon\R^{N}\backslash\{0\}\rightarrow (0,+\infty)$ satisfy $(H_1)$--$(H_3)$, where $s\in(0,1)$.  Assume also that $N>2s$.  Let $\Omega\subset\R^{N}$ be a bounded domain with Lipschitz boundary. Then, there exists a positive constant $C=C(N,s,\theta,\Omega)$, where $\theta$ is given by $(H_2)$, such that
\begin{equation}
    \|u\|_{L^{r}(\Omega)}\leqslant C\|u\|_{X},
    		\quad\mbox{ for } u\in X_0,
    \label{X0Sobineq}
\end{equation}
for any $r\in [1,2^{*}]$, where $2^{*} = \dfrac{2N}{N-2s}.$ 
\end{proposition}
For future reference in this article,  we state the estimate used to obtain the inequality in (\ref{X0Sobineq}):
\begin{equation}\label{ineq3}
    \|u\|_{L^{2^{*}}(\Omega)}^2\leqslant
        C
        \int_{\mathbb{R}^{2N}} \frac{(u(x)-u(y))^2}{|x-y|^{N+2s}} dx\, dy,
            \quad\mbox{ for } u\in X_0,
\end{equation}
for some positive constant $C$ (see 
\cite[Lemma 6]{Ser3}).

In view of Proposition \ref{X0EmbedThm}, the space $X_0$ is continuously embedded in
$L^r(\Omega)$ for any $r\in[1,2^{*}]$.  In fact, it follows from (\ref{ineq3}) and $(H_2)$ that 
\begin{equation}
    \|u\|_{L^{2^*}(\Omega)}^{2}\leqslant \frac{C}{\theta} \int_{\R^{2N}}(u(x)-u(y))^2K(x-y)\,dx\,dy
    \nonumber 
\end{equation}
for $u\in X_{0};$ so that, 
\begin{equation}
    \|u\|_{L^{2^*}(\Omega)}^2\leqslant \frac{C}{\theta}\|u\|_{X_0}^2,\quad\mbox{ for }u\in X_{0},
    \label{eineq1}
\end{equation}
where $\|u\|_{X_0}$ is the norm in $X_0$ defined in (\ref{normx0}). 

We can also deduce from (\ref{eineq1}) in which $\Omega$ is bounded with Lipschitz boundary that 
\begin{equation}
    \|u\|_{L^r(\Omega)}\leqslant C(\theta)\|u\|_{X_0},\quad\mbox{ for }u\in X_{0},
    \label{eineq2}
\end{equation}
for all $r\in [1,2^*].$

The embedding
$X_0\subset L^r(\Omega)$ can be shown to be compact
for $r\in[1,2^{*})$. For a proof of this property and
more  properties of the space $X_0$, we refer to \cite[Sections $1.5.1$ and $1.5.2$]{Ser2}, or \cite{Ser3}. 
For more information on fractional Sobolev spaces, we refer to Di Nezza {\it et al.}\cite[page $524$]{Nezza1}.

In \cite{Ser3}, Servadei and Valdinoci studied the 
problem of existence of nontrivial solutions of the the following boundary value problem 
\begin{equation}\label{prob_1}
   \left\{\begin{aligned} -\mathcal{L}_{K} u(x)-\lambda u(x) & =f(x,u),  \text{ for } x\in\Omega; \\ 
    u & =0 \quad \text{ in } \R^{N}\backslash\Omega,
    \end{aligned} \right.
\end{equation}
where $\lambda\in\R$ and $\mathcal{L}_{K}$ is the 
integro-differential operator defined by 
\begin{equation}
    \mathcal{L}_{K}u(x) = \int_{\R^{N}}\left[(u(x+y)+u(x-y)-2u(x)\right]K(y)\,dy,\quad\mbox{ for }x\in\R^{N},
    \label{defop}
\end{equation}
and for $u\in X$, 
with $K:\R^{N}\backslash\{0\}\rightarrow (0,\infty)$ being a function that satisfies the assumptions $(H_1)$--$(H_3)$ and $f$ satisfies the following conditions:
\begin{itemize}
    \item [$(H_4)$] $f\in C^{1}(\overline{\Omega}\times\R,\R)$.
    \item [$(H_5)$] There exist $a_1,a_2>0$ and $q\in (2,2^{*}),\, 2^{*}:=2N/(N-2s)$, such that 
    $$|f(x,t)|\leqslant a_1+a_2|t|^{q-1},\quad\mbox{ for any }x\in\Omega, \, t\in\R.$$
    \item [$(H_6)$] $\displaystyle\lim_{t\rightarrow 0}\frac{f(x,t)}{t}=0,$ uniformly in $x\in\Omega.$
    \item [$(H_7)$] $tf(x,t)\geqslant 0$ for any $x\in\Omega$, $t\in\R.$
    \item [$(H_8)$] There exist $\mu>2$ and $r>0$ such that 
    $$0 < \mu F(x,t) \leqslant tf(x,t),$$
    for any $x\in\overline{\Omega}$ and $t\in\R$ such that $|t|\geqslant r$,
    where $F(x,t)=\int_{0}^{t}f(x,\tau)\,d\tau$ is the primitive of $f$.
\end{itemize}

We note that the case in which 
 $$K(x)=\theta |x|^{-(N+2s)},\quad\mbox{ for }x\in\R^{N}\backslash\{0\},$$
 with $\theta$ given by 
 \begin{equation}\label{ARTheta}
     \theta = \frac{2^{2s} \Gamma((N+2s)/2)}{2\pi^{N/2}|\Gamma(-s)|},
 \end{equation}
 where 
 $\Gamma$ denotes the gamma function\footnote{$\displaystyle \Gamma(z) = \int_0^\infty t^{z-1} e^{-t}\ dt,\ $ for $\mbox{Re}(z)>0$, and for 
 complex values other than $0$ or the negative integers $\Gamma$ is the analytic continuation of the integral formula.}, corresponds to the fractional Laplacian
 $ \mathcal{L}_{K} = -(-\Delta )^s$. See Kwa\'snicki \cite{kwa1} for a detailed approach to the definition of the fractional Laplacian. 

The spectrum of the operator $-\mathcal{L}_{K}$ was studied in \cite{Ser3}.  We list here some of
the properties of the spectrum of
$-\mathcal{L}_{K}$
that will be useful in this work.

Denote by $(\lambda_k)_{k\in\N}$ the sequence of eigenvalues 
$$0 < \lambda_{1} < \lambda_{2}\leqslant\lambda_{3}\leqslant\ldots\leqslant \lambda_{k}\ldots$$ 
of the problem 
\begin{equation}\label{prob_2}
   \left\{\begin{aligned} -\mathcal{L}_{K} u & =\lambda u; \quad \text{in } \Omega, \\ 
    u & =0, \quad\quad \text{in } \R^{N}\backslash\Omega,
    \end{aligned} \right.
\end{equation}
and $\varphi_k$ is an eigenfunction associated with $\lambda_k$, for $k=1,2,\ldots.$ 

It was shown in \cite[Proposition $9$]{Ser3} that $\lambda_k\rightarrow\infty$ as $k\rightarrow\infty$. It was also shown in 
\cite{Ser3} that the first eigenvalue of problem (\ref{prob_2}) is simple and can be characterized by 
\begin{equation}
    \lambda_1 = \min_{u\in X_{0}\backslash\{0\}}\frac{\displaystyle\int_{\R^{2N}}|u(x)-u(y)|^2K(x-y)\,dx\,dy}{\displaystyle\int_{\Omega}|u(x)|^2\,dx}.
    \label{firsteigcharact}
\end{equation}

The characterization of the first eigenvalue $\lambda_1$ of the operator $-\mathcal{L}_{k}$ in (\ref{firsteigcharact}) implies that  
\begin{equation}
    \lambda_1\|u\|_2^2  \leqslant \|u\|_{X_0}^2,\quad\mbox{ for }u\in X_0.
    \label{pp1}
\end{equation}

The main result proved in \cite{Ser3} is the following theorem.

\begin{theorem}\cite[Theorem $1$]{Ser3}
    Let $s\in (0,1), N>2s$ and $\Omega$ be an open bounded set of $\R^{N}$ with Lipschitz boundary. Let $K:\R^{N}\backslash\{0\}\rightarrow (0,+\infty)$ be a function satisfying conditions $(H_1)$--$(H_3)$ and let $f:\overline{\Omega}\times\R\rightarrow\R$ verify $(H_4)$--$(H_8)$. Then, for any $\lambda\in\R$, problem (\ref{prob_1}) admits a solution $u\in X_{0}$ that is not identically zero. 
    \label{sertheo1}
\end{theorem}

By a solution of problem (\ref{prob_1}) we mean a 
weak solution; that is, a function $u\in X_0$ 
satisfying 
 \begin{equation} 
 \begin{aligned}
\int_{\R^{2N}}\left(u(x)-u(y)\right)\left(\varphi(x)-\varphi(y)\right)K(x-y)\,dx\,dy - \lambda\int_{\Omega}u(x)\varphi(x)\,dx = \\   \int_{\Omega}f(x,u(x))\varphi(x)\,dx,
\end{aligned}
     \label{weakform}
 \end{equation}
for all $\varphi\in X_{0}$.

In this work, our goal is to extend the multiplicity results in Theorem \ref{sertheo1}. The first multiplicity result obtained in this article deals with the case $ \lambda<\lambda_1$ in problem (\ref{prob_1}), where $\lambda_1$ is the first eigenvalue of the operator $-\mathcal{L}_{K}$. In \cite{Ser3}, Servadei and Valdinoci proved the existence of two nontrivial solutions given
by the mountain-pass theorem of Ambrosetti and Rabinowitz \cite{AmbRab}, where one is non-positive and the other is non-negative. Using infinite-dimensional Morse theory, we extend that result in the following theorem. 

\begin{theorem}
    Let all the assumptions of Theorem \ref{sertheo1} be satisfied. Then, for any $\lambda < \lambda_1$, problem (\ref{prob_1}) admits at least three nontrivial solutions where one is non-positive and another is non-negative. 
    \label{maintheoint}
\end{theorem}

For the case $\lambda\geqslant\lambda_1$ in problem (\ref{prob_1}), the authors of \cite{Ser3} proved the existence of a nontrivial solution $u_o$ of problem (\ref{prob_1}) by showing that the associated energy functional satisfies the hypotheses of the linking theorem of Rabinowitz \cite{Rab78}. We will extend this result by computing the critical groups of the associated energy functional at the origin and showing that the origin satisfies a local linking condition with respect to a particular decomposition of a Hilbert space $X_0=X_0^{-}\oplus (X_{0}^{-})^{\perp},$ where $k_o:=\dim X_{0}^{-}<\infty$. To obtain this result,
we need to replace hypothesis ($H_5$) by the following
growth condition on $f_t(x,t)$:
\begin{itemize}
    \item[($H_{9}$)] We assume that there exist $c_1,c_2>0$ such that 
    \begin{equation} 
    |f_t(x,t)|\leqslant c_1+c_2|t|^{q-2},
    \quad\mbox{ for all } x\in\Omega \mbox{ and }
    t\in\R.
\end{equation}
\end{itemize}

We will use the following decomposition of the Hilbert space $X_0 = X_{0}^{-}\oplus X_{0}^{+}$ with 
\begin{equation}
    X_{0}^{-}=\bigoplus_{j=1}^{k} \ker(-\mathcal{L}_{K}-\lambda_j I) \quad \mbox{ and }\quad X_{0}^{+}=(X_{0}^{-})^{\perp},
    \label{hilbdecomp}
\end{equation}
where $\lambda_{k}$  is the $k^{\mbox{th}}$ eigenvalue of the nonlocal operator $-\mathcal{L}_K$, for $k=1,2,\ldots$. For more details on the existence of the eigenvalues of the operator $-\mathcal{L}_K$ we refer to \cite[Proposition $9$]{Ser3}.

For $\lambda\geqslant\lambda_1$, we get the following multiplicity result. 

\begin{theorem}\label{maintheopaper22}
     Let $s\in (0,1), N>2s$ and $\Omega$ be an open bounded set of $\R^{N}$ with Lipschitz boundary. Let $K:\R^{N}\backslash\{0\}\rightarrow (0,+\infty)$ be a function satisfying conditions 
     $(H_1)$--$(H_3)$.  Assume that $f$ satisfies  hypotheses $(H_4)$, $(H_6)$--$(H_9)$.  Let $u_o$ be the solution found in Theorem \ref{sertheo1} through the linking theorem. Assume also that  $\lambda\in [\lambda_{k},\lambda_{k+1})$, for some  $k\in\N$, and define $X_0^{-}$ as in  (\ref{hilbdecomp}).  Then,      
     problem (\ref{prob_1}) admits at least two nontrivial solutions provided that
         \begin{equation}
         k_o+1\not\in[\mu_o,\mu_o+\nu_o],
         \nonumber
         \end{equation}
         where $k_o=\dim X_{0}^{-}$, $\mu_o$ and $\nu_0$ are the Morse index and nullity of $u_o$, respectively. 
\end{theorem}

The last result of this paper was motivated by the work of Servadei in \cite{Ser4}. In that article, the author proved the existence of infinitely many weak solutions to the problem (\ref{prob_1}) with the hypothesis $(H_8)$ replaced by the following: 
\begin{itemize}
\item[]  Assume that $f(x,t)$ is 
differentiable with respect to the second variable
and that there exist $\beta\in(0,1)$ and $r>0$ such that 
\begin{equation}\label{Servadei1.14}
    0<\frac{f(x,t)}{t}\leqslant \beta f_t(x,t),
\end{equation}
for any $x\in\overline{\Omega}$ and  $t\in\R$ such that $|t|\geqslant r$,
where $f_t$ denotes the derivative of $f$ with respect to the second variable.
\end{itemize}
In this work, we analyze the case in which (\ref{Servadei1.14}) is replaced by hypothesis $(H_8)$, and we assume that the nonlinearity $f$ is odd in the second variable. Then, applying the $\Z_2$ version of the mountain-pass theorem of Ambrosetti and Rabinowitz \cite{AmbRab}, we prove the following multiplicity result:

\begin{theorem}
     Suppose that 
  $K:\R^{N}\backslash\{0\}\rightarrow (0,+\infty)$  satisfies conditions $(H_1)$--$(H_3)$, and $f$ satisfies $(H_4),$ $(H_5),(H_8)$ and $f(x,s)$ is odd in s. Assume also that $\lambda\in [\lambda_{k},\lambda_{k+1})$, for some $k\in\N$, and $X=V\oplus W$ with $V=X_{0}^{-}$ and $W=V^{\perp}$ as in (\ref{hilbdecomp}). Then, problem (\ref{prob_1}) possesses an unbounded sequence of weak solutions.
    \label{secexistheo33}
\end{theorem}

A particular case of problem (\ref{prob_1}) for $s\in (0,1)$ is given by the kernel $K(x)=\theta|x|^{-(N+2s)}$, $x\in\R^{N}\backslash\{0\}$, 
where $\theta$ is given in (\ref{ARTheta}).  In this case, the operator $-\mathcal{L}_{K}$ becomes the fractional Laplacian $-\mathcal{L}_K = -(-\Delta)^s$.  Problems involving fractional Laplacian have been the subject of research for many years.  We refer the reader to the works of Di Nezza {\it et al.} in \cite{Nezza1}, Valdinoci in \cite{Vald1}, Bisci {\it et al.} in \cite{Ser2}, Caffarelli in \cite{Caf1}, and references therein, for a series of applications of this kind of operator and challenges faced when working with them.

One of the main motivations to extend the results in Theorem \ref{sertheo1} was the technique used by Wang in \cite[Section $3$]{Wang}. In that article, Wang extended the results produced by Ambrosetti and Rabinowitz in \cite{AmbRab} for the semilinear elliptic boundary value problem 
\begin{equation}\label{prob_w}
   \left\{\begin{aligned} -\Delta u & =f(u); \quad \text{in } \Omega, \\ 
    u & =0, \quad\quad \text{on } \partial\Omega,
    \end{aligned} \right.
\end{equation}
under similar assumptions to those in 
$(H_4)$--$(H_8)$. Using infinite-dimensional Morse theory, Wang extended the results of \cite{AmbRab} by proving the existence of at least three nontrivial solutions of problem (\ref{prob_w}). The argument relied on the computation of the critical groups of the associated energy functional for the problem (\ref{prob_w}) at the origin and at each of its mountain pass critical points. The existence of a 
third nontrivial solution was established using an argument by contradiction involving a long exact sequence of singular homology groups of a particular topological pair. Motivated by this result, we will extend the results of Theorem \ref{sertheo1} to the case in which $\lambda < \lambda_1$. 

For the case $\lambda\geqslant\lambda_1$, we will analyze the critical point $u_o$ found through the linking theorem of Rabinowitz in Theorem \ref{sertheo1}.  We consider the cases in which $u_o$ is nondegenerate and $u_o$ is degenerate separately. A result due to Gromoll and Meyer in \cite{Gr}, the shifting theorem, and the generalized Morse lemma will play a key role in proving the multiplicity result in the case in which $u_o$ is a degenerate critical point.

This paper is organized as follows: Section \ref{prelsection} presents some notation and preliminary results that will be used throughout this work. In Section \ref{coriginsection}, we compute the critical groups of the energy functional  associated with problem (\ref{prob_1}) at the origin for the two
cases: $\lambda < \lambda_1$ and $\lambda\geqslant\lambda_1.$ In the first case, we prove that the origin is a local minimizer for $J_\lambda$, the energy functional associated with problem (\ref{prob_1}), and for the second case we show that $J_\lambda$ satisfies a local linking condition at the origin introduced by Li and Liu in \cite{LiLiu2}. The computations of the critical groups of $J_\lambda$ at infinity are done in Section \ref{cinfitysection}. In Section \ref{secsolsection}, we prove the multiplicity results for problem (\ref{prob_1}).  In Section \ref{secexistsec}, we prove that the associated energy functional for problem (\ref{prob_1}) has an unbounded sequence of critical points when $\lambda\geqslant\lambda_1 $ and $f(x,s)$ is odd in the second variable $s$. The result is obtained through an application of the $\Z_2$ version of the 
mountain-pass theorem of Ambrosetti and Rabinowitz  in \cite{AmbRab}.

\section{Preliminaries}\label{prelsection}

In this section, we present some preliminary concepts and estimates that will be used throughout this work.

We will employ a variational approach to analyze problem (\ref{prob_1}). Let $X_0$ be the Hilbert space defined in (\ref{x0def}) and   $J_\lambda:X_0\rightarrow\R$ be the energy functional associated with problem\eqref{prob_1} given by
\begin{equation}
    J_\lambda(u):=\frac{1}{2}\|u\|_{X_0}^2 -\frac{\lambda}{2} \int_{\Omega} u^2 dx-\int_\Omega F(x,u)dx, \quad\mbox{ for }u\in X_{0},
    \label{func0}
\end{equation}
where $\|\cdot\|_{X_0}$ is the norm defined in (\ref{normx0}). 

The functional $J_\lambda$ is well-defined by virtue of hypotheses $(H_4)$--$(H_5)$ and the embedding theorem in Proposition \ref{X0EmbedThm}. $J_\lambda$ is also Fr\'{e}chet differentiable at every $u\in X_{0}$ and its 
Fr\'{e}chet derivative at $u$ is given by 
\begin{equation}
\begin{array}{rcl}
\langle J_{\lambda}^{\prime}(u),\varphi\rangle
& = & \displaystyle
\int_{\R^{2N}}(u(x)-u(y))(\varphi(x)-\varphi(y))K(x-y)\,dx\,dy\\
					\\
&   & \displaystyle \qquad
-\ \lambda\int_{\Omega}u\varphi\,dx-\int_{\Omega}f(x,u)\varphi\,dx,
\end{array}
\label{frechet0}
\end{equation}
for all $\varphi\in X_{0}$, where $\langle\cdot,\cdot\rangle:X_{0}^{*}\times X_{0}\rightarrow\R$ denotes the duality pairing and $X_{0}^{*}$ is the dual space of $X_0.$

Next, we present some estimates for the functions $f$ and its primitive $F$. It follows from  \cite[Lemma $3$]{Ser3} that the assumptions $(H_4)$--$(H_7)$ imply that, for any $\varepsilon >0$, there exists $\delta=\delta(\varepsilon)>0$ such that, for any $x\in\Omega$ and $t\in\R,$ we have
\begin{equation}
    |f(x,t)|\leqslant 2\varepsilon|t| + q\delta(\varepsilon)|t|^{q-1};\quad\mbox{ for }(x,t)\in \Omega\times\R;
    \label{estlittlef}
\end{equation}
so that, integrating (\ref{estlittlef}),
\begin{equation}
    |F(x,t)|\leqslant \varepsilon |t|^2 + \delta(\varepsilon)|t|^{q},\quad\mbox{ for }(x,t)\in \Omega\times\R.
    \label{estbigF}
\end{equation}

The hypotheses $(H_4)$ and $(H_8)$ imply the existence of constants $a_1,a_2 > 0$ such that
\begin{equation}
    F(x,t)\geqslant a_{1}|t|^{\mu}-a_2,
    \quad\mbox{ for } x\in\overline{\Omega}
        \mbox{ and } t\in\R.
    \label{estARbigF}
\end{equation}
For a detailed proof of (\ref{estARbigF}), we refer to \cite[Lemma $4$]{Ser3}. The estimate (\ref{estARbigF}) shows that $F(x,t)$ has a superquadratic growth in $t$ for any $x\in\overline{\Omega}$ and $t\in\R$.  Nonlinearities with growth of this type were first considered by Ambrosetti and Rabinowitz in \cite{AmbRab}. The estimate in (\ref{estARbigF})  is needed in the proof
of the fact that the functional $J_\lambda$ given in
(\ref{func0}) satisfies the Palais-Smale condition,
which was first introduced in \cite{SM2} and presented here as found in Perera and Schechter \cite{PerSch}.

\begin{definition}\cite[Page $2$]{PerSch}
   {\rm  A functional $J\in C^{1}(X_0,\R)$ satisfies the Palais-Smale compactness condition at the level $c$, or $PS_c$ for short, if every sequence $(u_j)\subset X_0$ such that 
    $$J(u_j)\rightarrow c\quad\mbox{ and }\quad J^{\prime}(u_j)\rightarrow 0\quad\mbox{ as }j\rightarrow\infty,$$
called a $PS_c$ sequence, has a convergent subsequence. $J$ satisfies $PS$ if it satisfies $PS_c$ for every $c\in\R$, or equivalently, if every sequence such that $J(u_j)$ is bounded and $J^\prime(u_j)\rightarrow 0$ as $j\rightarrow\infty$, called a PS sequence, has a convergent subsequence.} 
\end{definition}

It was shown in \cite{Ser3} that, under the 
assumptions $(H_1)$, $(H_2)$, $(H_3)$, $(H_4)$ and
$(H_8)$, the functional $J_\lambda$ given in 
(\ref{func0}) satisfies the Palais-Smale condition. Servadei and Valdinoci proved this result by analyzing the two  cases: $\lambda < \lambda_{1}$ \cite[Propositions $13$, $14$]{Ser3} and  $\lambda \geqslant \lambda_1$ separately \cite[Propositions $14$, $20$]{Ser3}. 

The PS condition is also needed to apply some of the results from infinite-dimensional Morse theory. Before presenting those results, we first introduce the concept of critical groups.

Define $J_\lambda^{c}=\{w\in X_0|\, J_\lambda(w)\leqslant c\}$,  the sub--level set of $J_\lambda$ at $c$, and set 
$$
    \mathcal{K}=\{u\in X_{0}|\, J_\lambda^{\prime}(u)=0\} ,
$$ 
the critical set of $J_\lambda$.  For  an isolated critical point  $u_0$  of $J_\lambda$, the $k^{\mbox{th}}$ critical groups of $J_\lambda$ at $u_0$, with coefficients in a field $\F$, are defined by 
\begin{equation}
    C_{k}(J_\lambda,u_0)=H_{k}(J_\lambda^{c_{0}}\cap U,J_\lambda^{c_{0}}\cap U\backslash\{u_0\}),\quad\mbox{ for all }k\in\Z,
    \label{cgdef}
\end{equation}
where $c_{0}=J_\lambda(u_0)$, $U$ is a neighborhood of $u_0$ that contains no critical points of $J_\lambda$ other than $u_0$, and $H_{*}$ denotes the singular homology groups. The critical groups are independent of the choice of $U$ by the excision property of homology (see  Hatcher \cite[Theorem 2.20]{AH}). For more information on the definition of critical groups, we refer the reader to \cite{PerSch}, \cite{KC},  \cite{MMP}, and \cite{MW1}. 

We will also refer to the reduced singular homology groups  $\widetilde{H}_{*}$ in the final sections of this work. The reduced homology groups of a particular topological space $Y$ are related to its homology groups by 
\begin{equation}
    H_{0}(Y) = \widetilde{H}_{0}(Y)\oplus \F\quad\mbox{ and }\quad \widetilde{H}_{k}(Y)\cong H_{k}(Y),\mbox{ for }k>0.
    \label{redhomrel}
\end{equation}
For more details, see \cite[Page $110$]{AH}.

The solutions of problem (\ref{prob_1}) given in 
Theorem \ref{sertheo1} were obtained as critical points of $J_\lambda$ via the mountain pass theorem of Ambrosetti and Rabinowitz \cite{AmbRab}.

Next, we present the definition of a mountain-pass type critical point given by Hofer in \cite{Hofer2} 
and found in Montreanu {\it et al.} \cite{MMP}.

\begin{definition}\cite[Definition $6.98$]{MMP} Let $E$ be a Banach space, $J\in C^{1}(E,\R)$, and $u_0\in\mathcal{K}.$ We say that $u_0$ is of {\it mountain-pass type} if, for any open neighborhood $U$ of $u_0$, the set $\{w\in U|\, J(w)<J(u_0)\}$ is nonempty and not path-connected.    
\label{mtpassdef}
\end{definition}

The following is a variant of the mountain pass theorem of Ambrosetti and Rabinowitz due to Hofer \cite{Hofer2}. We present the result here as stated in Montreanu {\it et al.} \cite{MMP}. 

\begin{theorem}\cite[Theorem $6.99$]{MMP}
    Assume that $E$ is a Banach space and that $J\in C^{1}(E,\R)$ satisfies the Palais-Smale condition. 
 Let $u_o, u_1\in E$ and define 
 $$\Gamma=\{\gamma\in C([0,1],E)\,|\,\gamma(0)=u_o,\gamma(1)=u_1\}.$$
 Put
$$c:=\displaystyle\inf_{\gamma\in\Gamma}\max_{t\in[0,1]}J(\gamma(t)).$$ 
If $c>\max\{J(u_o),J(u_1)\}$, then $\mathcal{K}^{c}\ne\emptyset$. 

Moreover, 
if $\mathcal{K}^c=\{u\in \mathcal{K} \,|\, J(u)=c\}$ is discrete, there exists $u\in\mathcal{K}^c$ that is of mountain-pass type.
    \label{hoftheo2}
 \end{theorem}

In \cite{Ser3}, Servadei and Valdinoci proved the existence of two critical points for the functional $J_\lambda$ through the mountain-pass theorem of Ambrosetti and Rabinowitz when $\lambda < \lambda_1$. It can be shown that the hypotheses of Theorem \ref{sertheo1} also imply the hypotheses of Theorem \ref{hoftheo2}.  In Section \ref{secexistsec}, we will use an argument by contradiction in which we assume the set $\mathcal{K}$ consists of only three critical points; so that $\mathcal{K}$ is discrete. Then, at least one of the  critical points found by Servadei and Valdinoci in \cite{Ser3} is of mountain-pass type per Theorem \ref{hoftheo2}. This will allow us to obtain information about the $1^{\mbox{st}}$ critical group of $J_\lambda$ at a mountain-pass type critical point as shown in the next proposition.

\begin{proposition}
    \cite[Proposition $6.100$]{MMP} Let $E$ be a reflexive Banach space, $J\in C^{1}(E,\R),$ and $u_0\in\mathcal{K}$ be an isolated critical point of $J$. If $u_0$ is of mountain-pass type, then $C_{1}(J,u_0)\not\cong 0.$
    \label{propMMP}
\end{proposition}

We will prove the existence of a second nontrivial solution of problem (\ref{prob_1}) using a result due to Chang \cite{KCM} as presented in Perera and Schechter \cite{PerSch}. 

\begin{proposition}\cite[Proposition $1.4.1$]{PerSch}
Let $a$ and $b$ be extended real numbers satisfying 
$-\infty < a < b\leqslant +\infty$.  Suppose that $J\in C^1(E,\R)$ has a finite number of critical points at the level $c\in (a,b)$, and that there are no other critical values of $J$ in $[a,b]$.  Assume also that $J$ satisfies the $PS_{c^\prime}$ condition for all $c^\prime\in [a,b]$. Then, 
\begin{equation}
    H_{k}(J^b,J^{a})\cong \bigoplus_{u\in \mathcal{K}^{c}}C_{k}(J,u),\quad\mbox{ for all }k\in \Z,
    \label{perprop14}
\end{equation}
where $\mathcal{K}^c=\{u\in \mathcal{K}\, |\, J(u)=c\}.$ 
\label{mainprop}
\end{proposition}

Finally, we present some concepts and results that will be used in the proof of the multiplicity result for problem (\ref{prob_1}) in the case $\lambda\geqslant\lambda_1.$
 
\begin{definition}\cite[Chapter 5, page $101$]{SchecBookFunc}
    {\rm
     Let $E,Y$ be Banach spaces, $A:E\rightarrow Y$ be a bounded linear operator,  and $A^{\prime}$ be its adjoint operator. $A$ is said to be a {\it Fredholm operator} from $E$ to $Y$ if 
     \begin{itemize}
         \item [(a)] $\dim N(A)$ is finite,
         \item [(b)] $R(A)$ is closed in $Y$,
         \item [(c)] $\dim N(A^{\prime})$ is finite,
     \end{itemize}
    }
    where $N(A)$ denotes the kernel of $A$ and $R(A)$ denotes the range of $A$.
    \label{deffredd}
\end{definition}

Let $X_0$ be the Hilbert space defined in (\ref{x0def}), $u\in X_0$ and $J \in C^{2}(X_0,\R)$. We can define implicitly the linear operator $L(u):X_0\rightarrow X_0$ by
\begin{equation} J^{\prime\prime}(u)(v,w)=\langle L(u)v,w\rangle _{X_0};\quad\mbox{ for }v,w \in X_0,
\label{isoope}
\end{equation}
by virtue of the Riesz-Fr\'{e}chet representation theorem (see Br\'{e}zis \cite[Theorem $5.5$]{Brezis}). It can be shown that $L(u)$ is a self-adjoint operator and we can identify $L(u)$ with $J^{\prime\prime}(u).$  In the case where $J^{\prime\prime}(u)$ is a Fredholm operator, the Hilbert space $X_0$ is the orthogonal sum of $R(J^{\prime\prime}(u))$ and $\ker J^{\prime\prime}(u).$ For more details, we refer to \cite[Pages 184, 185]{MW1}.

\begin{definition}
{\rm 
Assume that $J\in C^{2}(X_0,\R)$ and $u_o$ is a critical point of $J$. The {\it Morse index} of $u_o$, denoted by $\mu_o$, is defined as the supremum of the dimensions of the vector spaces of $X_0$ on which $J^{\prime\prime}(u_o)$ is negative definite. The {\it nullity} of $u_o$, denoted by $\nu_o$, is defined as the dimension of $N( J^{\prime\prime}(u_o)).$ Set $\mu^{*}=\mu_o+\nu_o$; 
this is called the {\it large Morse index} of
$u_o$.
}
\label{defmorseindex}
\end{definition}

\begin{definition}
  {\rm   Let $E$ be a Hilbert space and $E^{*}$ be its dual. A critical point $u_o$ of $J\in C^{2}(E,\R)$ is called {\it nondegenerate} if $J^{\prime\prime}(u_o):E\rightarrow E^{*}$ is an isomorphism. Otherwise, $u_o$ is called a {\it degenerate }critical point of $J$. 
  }
\end{definition}
 
The following result is a consequence of a result due to Gromoll and Meyer in \cite{Gr} and the generalized Morse lemma, also known as the splitting lemma (see \cite[Theorem $8.3$] {MW1}).  We present the result as found in Cingolani and Vannella \cite{Cing2}.

\begin{theorem}\cite[Proposition $2.5$]{Cing2}
Suppose that $E$ is a Hilbert space and $J\in C^{2}(E,\R)$. Let $u_o$ be an isolated critical point of $J$ with Morse index $\mu_o$ and large Morse index $\mu^\ast.$ Suppose that $J^{\prime\prime}(u_o)$ is a Fredholm operator and let $V$ be the kernel of $J^{\prime\prime}(u_o)$. 
\begin{itemize}
    \item [(a)] If $u_o$ is a local minimum of $J|_{V}$ then 
    $$C_{k}(J,u_o)\cong \delta_{k,\mu_o}\F,\quad\mbox{ for }k\in\Z.$$

    \item [(b)] If $u_o$ is a local maximum of $J|_{V}$ then 
      $$C_{k}(J,u_o)\cong \delta_{k,\mu^{\ast}}\F,\quad\mbox{ for }k\in\Z.$$

    \item [(c)] If $u_o$ is neither a local maximum nor a local minimum of $J|_{V}$, then 
    $$C_{k}(J,u_o)\cong 0,\quad\mbox{ for }k\not\in (\mu_0,\mu^\ast).$$
\end{itemize}
    \label{cingotheo}
\end{theorem}

The previous theorem will be applied in the proof of Theorem \ref{maintheopaper} to obtain a multiplicity result for problem (\ref{prob_1}) when $\lambda\geqslant \lambda_1$.

\begin{remark} \label{fredremark}
{\rm 
Let $u_o$ be the critical point of $J_\lambda$ found in Theorem \ref{sertheo1} through the linking theorem of Rabinowitz for the case $\lambda\geqslant \lambda_1$. By virtue of the hypotheses $(H_4)$ and $(H_9)$, the second Fr\'{e}chet derivative $J_\lambda^{\prime\prime}(u_o)$ is a well-defined bilinear form in the Hilbert space $X_0$ that is given by 
\begin{equation}
    J_\lambda^{\prime\prime}(u_o)(v,w)=\langle v,w\rangle_{X_0} - \lambda \langle v,w\rangle _{L^2} - \int_{\Omega}f_t(x,u_o)vw\,dx,
    \label{secfrec}
\end{equation}
for all $v,w\in X_0$. 

We show that $J_\lambda^{\prime\prime}(u_o)$ is a Fredholm operator. In fact, it follows from (\ref{isoope}) that there exists a self-adjoint operator $L(u_o):X_0\rightarrow X_0$ such that we can write (\ref{secfrec}) as 
\begin{equation}
   \langle L(u_o)v,w\rangle_{X_0} =\langle v,w\rangle_{X_0} -  \langle \lambda v + f_t (\cdot,u_0)v,w\rangle _{L^2(\Omega)},
    \label{opdef0}
\end{equation}
for all $v,w\in X_0$, where 
\begin{equation}
    L(u_o)(v) = v - \mathcal{N}_{u_o}(v), \quad\mbox{ for }v\in X_{0}, 
    \label{opdef}
\end{equation}
where $\mathcal{N}_{u_o}$ is a compact operator in 
$X_0$ by the compactness of the embedding $X_0\hookrightarrow L^q(\Omega)$, since $q<2^\ast$.
This fact follows from Theorem \ref{X0EmbedThm}. 

Then, the operator $L(u_o)$, defined in (\ref{opdef}), is of the form $I-\mathcal{N}_{u_o}$ where $I:X_{0}\rightarrow X_0$ is the identity operator and $\mathcal{N}_{u_o}:X_0\rightarrow X_0$ a compact operator. Therefore, the operator $L(u_o)$ is a compact perturbation of the identity. Then, it follows from the Fredholm alternative theorem (see Schechter \cite[Theorem 4.12]{SchecBookFunc}) that $L(u_o)$ is a Fredholm operator from which we conclude that $J_\lambda^{\prime\prime}(u_o)$ is a Fredholm operator. 
 
}
\end{remark}

\section{Critical Groups at the Origin}\label{coriginsection}

In this section, we will compute the critical groups of $J_\lambda$ at the origin for the cases $\lambda < \lambda_1$ and $\lambda \in  [\lambda_k,\lambda_{k+1}),$ for some $k\geqslant 1,$ respectively. We have to treat both cases separately because of the different geometries of the level sets of the functional $J_\lambda$ at the origin
associated with the two cases.

\subsection{Case 1: } To compute the critical groups of $J_\lambda$ at the origin in the case $\lambda<\lambda_1$ we will use an auxiliary estimate proved by Servadei and Valdinoci in \cite{Ser3}.

\begin{lemma}\cite[Lemma 10]{Ser3}
Let $K:\R^{N}\backslash\{0\}\rightarrow (0,\infty)$ satisfy assumptions $(H_1)$--$(H_3)$ and let $\lambda < \lambda_1$. Then, there exist two positive constants $m_1^\lambda$ and $M_1^\lambda$, depending only on $\lambda,$ such that 
\begin{equation}
m_1^\lambda\|u\|_{X_0}^2\leqslant\int_{\R^{2N}}|u(x)-u(y)|^2K(x-y)\,dx\,dy - \lambda\int_{\Omega}|u(x)|^2\,dx \leqslant M_1^\lambda\|u\|_{X_0}^2,
    \label{ineq11}
\end{equation}
for any $u\in X_0$.
    
\end{lemma}

In the next lemma, we will compute the critical groups of $J_\lambda$ at the origin for the case 
$\lambda < \lambda_1$. 

\begin{lemma}\label{lemmacorigincase1}
    Assume that $K:\R^{N}\backslash\{0\}\rightarrow (0,\infty)$ satisfies assumptions $(H_1)$--$(H_3)$
    and that $f$ satisfies $(H_4)$--$(H_8)$. 
    Assume also that $\lambda < \lambda_1$. Then, the origin in $X_0$ is a local minimizer for $J_\lambda$. Moreover, the  critical groups of $J_\lambda$ at the origin are given by
    \begin{equation}
        C_{k}(J_\lambda,0)\cong\delta_{k,0}\F,\quad\mbox{ for }k\in\Z.
        \label{corigin}
    \end{equation}
\end{lemma}
\begin{proof}[Proof:] 
We will show that the origin is a local minimizer for $J_\lambda$. In fact, by virtue of the estimate (\ref{estbigF}), for any $\varepsilon>0$, there exists
$\delta(\varepsilon)>0$ such that 
\begin{equation}
    \int_{\Omega}F(x,u)\,dx \leqslant \varepsilon\int_{\Omega}|u|^2\,dx + \delta(\varepsilon)\int_{\Omega}|u|^{q}\,dx,
        \quad\mbox{ for } u\in X_0.
    \label{est0001}
\end{equation}

We note that $u\in L^q(\Omega)$ for every $u\in X_0$
by virtue of Proposition \ref{X0EmbedThm} (see also
the estimate in (\ref{X0Sobineq})).

Next, use the assumption $N>2s$, in conjunction with
H\"{o}lder's inequality, with coefficients $p_{1}=N/(N-2s)$ and $q_{1}=N/2s$ in (\ref{est0001}), to get 
$$
	\int_{\Omega}F(x,u)\,dx  \leqslant
	\varepsilon\left(\int_{\Omega}|u|^{2^{*}}\,dx\right)^{\frac{N-2s}{N}}|\Omega|^{\frac{2s}{N}}+\delta(\varepsilon)\int_{\Omega}|u|^{q}\,dx;
$$
so that,
\begin{equation}
    \int_{\Omega}F(x,u)\,dx  \leqslant 
    \varepsilon \|u\|_{L^{2^*}(\Omega)}^2|\Omega|^{\frac{2s}{N}} + \delta(\varepsilon)\|u\|_{L^q(\Omega)}^q,
    \label{est0002}
\end{equation}
where $|\Omega|$ denotes the Lebesgue measure of the bounded set $\Omega$.  

Next, since $q\in (2,2^*)$, we  use the estimate in (\ref{X0Sobineq}) to
obtain 
\begin{equation}
    \|u\|_{L^{q}(\Omega)}\leqslant C\|u\|_{X},
    		\quad\mbox{ for } u\in X_0
    \label{MyEst0005}
\end{equation}
and
\begin{equation}
    \|u\|_{L^{2^{*}}(\Omega)}\leqslant C\|u\|_{X},
    		\quad\mbox{ for } u\in X_0,
    \label{MyEst0010}
\end{equation}
where $C$ is a positive constant depending on $N$, $s$, $\Omega$ and $\theta$. 

Substituting the estimates in (\ref{MyEst0005}) and
(\ref{MyEst0010}) on the right-hand side of the 
estimate in (\ref{est0002}) yields the estimate
\begin{equation}
    \int_{\Omega}F(x,u)\,dx   
       \leqslant \left(C_{1}\varepsilon + C_2\delta(\varepsilon)\|u\|_{X_{0}}^{q-2}\right)\|u\|_{X_0}^2,
       \label{ineq66}
\end{equation}
for $u\in X_0$, where $C_1$ and $C_2$ are positive constants depending on  $N$, $s$, $\Omega$ and $\theta$. 

Next, set 
\begin{equation}\label{rhoDfn}
    \rho=\displaystyle\left(\frac{\varepsilon}{C_{2}\delta(\varepsilon)}\right)^{1/(q-2)}.
\end{equation}
Since  $q>2$ we then obtain from  (\ref{ineq66}) that 
\begin{equation}
    \|u\|_{X_0} < \rho \quad \Rightarrow \int_{\Omega}F(x,u)\,dx \leqslant C_3\varepsilon\|u\|_{X_0}^2,
    \label{ineq77}
\end{equation}
for some positive constant $C_3$. 

Hence, by virtue of the definition of the functional $J_\lambda$ in (\ref{func0}), and the estimates  (\ref{ineq11}) and (\ref{ineq77}), we get  
\begin{equation}
    \begin{aligned}
        J_{\lambda}(u) & = \frac{1}{2}\|u\|_{X_0}^2 -\frac{\lambda}{2} \int_{\Omega} u^2 dx-\int_\Omega F(x,u)dx \\ 
        & \geqslant \left[\frac{m_1^\lambda}{2} - C_{3}\varepsilon \right]\|u\|_{X_0}^2,  
    \end{aligned}
    \label{ineq5}
\end{equation}
for $u\in X_{0}$, where $C_{3}$ is a positive constant. 

Choose 
$\varepsilon 
= \displaystyle \frac{m_1^\lambda}{4C_3} $ in \eqref{ineq5} and
$\rho$ as given in (\ref{rhoDfn}) to obtain 
\begin{equation}
    J_\lambda (u) \geqslant \frac{m_1^{\lambda}}{4}\|u\|_{X_0}^2 > J_\lambda(0),\quad\mbox{ for }0<\|u\|_{X_0} < \rho.
    \label{ineq88}
\end{equation}
Consequently, the origin is a local minimum of $J_\lambda$ in $B_{\rho}(0)$.  Therefore, by virtue of the definition of the critical groups of $J_\lambda$ at the origin in (\ref{cgdef}) with $U=B_{\rho}(0)$, it follows from \cite[Example $1$, page $33$]{KC} that 
$$C_{k}(J_\lambda,0)\cong \delta_{k,0}\F,\quad\mbox{ for }k\in \Z.$$
This concludes the proof of the lemma. 
\end{proof}

\subsection{Case 2:} 

We start with the notion of local linking introduced by Li and Liu in \cite{LiLiu2}, which will allow us to compute the critical groups of $J_\lambda$ at the origin for $\lambda\geqslant\lambda_1.$ We present the definition as found in Li and Willem \cite{LiWi} for the reader's convenience. 

\begin{definition}\cite[Section $0$]{LiWi}
    Let $J$ be a $C^1$ functional defined on a Banach space $X$. We say that $J$ has a local linking near the origin if $X$ has a direct sum decomposition $X=X^{-}\oplus X^{+}$ with $\dim X^{-}<\infty$, $J(0)=0$, and, for some $\rho >0$, 
    \begin{equation}
        \begin{aligned} 
        J(u) & \leqslant 0, \quad\mbox{ for }u\in X^{-}, \quad \|u\|\leqslant \rho, \\
        J(u) & > 0, \quad \mbox{ for }u\in X^{+},\quad  0 < \|u\| \leqslant \rho. 
        \end{aligned}
        \label{locallinkcond}
    \end{equation}
    \label{locallinkdef}
\end{definition}

We will need an auxiliary lemma to prove that $J_\lambda$ has a local linking at the origin with respect to the decomposition in (\ref{hilbdecomp}). 

\begin{lemma}\cite[Lemma $16$]{Ser3}
Let $K:\R^{N}\backslash\{0\}\rightarrow (0,\infty)$ satisfy the assumptions $(H_1)$--$(H_3)$ and let $\lambda\in[\lambda_k,\lambda_{k+1})$ for some $k\in\N$. Then, for any $v\in X_{0}^{+}$, 
\begin{equation}
    \int_{\R^{2N}}|v(x)-v(y)|^2 K(x-y)\,dx\,dy - \lambda \int_{\Omega}|v(x)|^2\,dx\geqslant m_{k+1}^{\lambda}\|v\|_{X_0}^2,
    \label{ineqspt}
\end{equation}
where $m_{k+1}^\lambda = 1 - \dfrac{\lambda}{\lambda_{k+1}} > 0.$
\label{serlema16}
\end{lemma}

Next, we show that the functional $J_\lambda$ has a local linking at the origin with respect to the decomposition $X_0=X_{0}^{-}\oplus X_{0}^{+},$ where $X_0^{-}$ and $X_{0}^{+}$ are defined in (\ref{hilbdecomp}). 

 \begin{lemma}
 Let $K:\R^{N}\backslash\{0\}\rightarrow (0,\infty)$ satisfy the assumptions $(H_1)$--$(H_3)$.
    Assume that there exists $k\geqslant 1$ such that $\lambda_k < \lambda_{k+1}$ and $\lambda\in [\lambda_k,\lambda_{k+1})$. In addition, assume that  $(H_4)$--$(H_8)$ are satisfied. Then, $J_\lambda$ has a local linking at the origin with respect to the decomposition $X_0=X_{0}^{-}\oplus X_{0}^{+}$ defined in (\ref{hilbdecomp}). Moreover, 
    \begin{equation}
        C_{k_o}(J_\lambda,0)\not\cong 0,
        \label{cgrouporigincase2}
    \end{equation}
    where $k_o=\dim X_{0}^{-}.$
    \label{lemmalocallink}
\end{lemma}

\begin{proof}[Proof:] We have to show that the functional $J_\lambda$ satisfies the conditions in (\ref{locallinkcond}). 

First, for $u\in X_{0}^{-}$, it was shown in \cite[Proposition $18$]{Ser3} that $J_\lambda(u)\leqslant 0$ for any $u\in X_{0}^{-}$. 

Next, we proceed to prove that $J_\lambda(u)>0$ for any $u\in X_{0}^{+}$ in a small neighborhood of the origin such that $u\not=0$. 

Using the definition of $J_\lambda$ in 
(\ref{func0}), and applying the estimates in 
(\ref{estbigF}) and (\ref{ineqspt}), as well as 
H\"{o}lder's inequality, we get 
\begin{eqnarray}
        J_\lambda (u)&\geqslant &\frac{1}{2}\|u\|_{X_0}^2-\frac{\lambda}{2}\|u\|_{L^2(\Omega)}^2-\varepsilon \|u\|_{L^2(\Omega)}^2 -\delta \|u\|_{L^q(\Omega)}^q \nonumber\\
    &\geqslant & \frac{m_{k+1}^\lambda}{2}\|u\|_{X_0}^2-\varepsilon \|u\|_{L^2(\Omega)}^2 -\delta \|u\|_{L^q(\Omega)}^q \nonumber\\
    &\geqslant & \frac{m_{k+1}^\lambda}{2}\|u\|_{X_0}^2-\varepsilon|\Omega|^\frac{2^*-2}{2^*} \|u\|_{L^{2^*}(\Omega)}^2 -\delta |\Omega|^\frac{2^*-q}{2^*}\|u\|_{L^{2^*}(\Omega)}^q 
\label{J>0_1}
\end{eqnarray}

Set 
\begin{equation}
    [u]^2_{s,2}:=\int_{\mathbb{R}^{2N}} \frac{(u(x)-u(y))^2}{|x-y|^{n+2s}} dx\, dy,
    \label{colcus2}
\end{equation}
and use (\ref{ineq3}) in (\ref{J>0_1}), as well as the notation defined in (\ref{colcus2}) , to obtain from (\ref{J>0_1}) that 
\begin{eqnarray}
      J_{\lambda}(u)  &\geqslant & \frac{m_{k+1}^\lambda}{2}\|u\|_{X_0}^2-\varepsilon C|\Omega|^\frac{2^*-2}{2^*} [u]_{s,2}^2 -\delta C |\Omega|^\frac{2^*-q}{2^*}[u]_{s,2}^q \nonumber\\
        &\geqslant & \frac{m_{k+1}^\lambda}{2}\|u\|_{X_0}^2-\frac{\varepsilon C}\theta|\Omega|^\frac{2^*-2}{2^*} \|u\|_{X_0}^2 -\frac{\delta C}\theta |\Omega|^\frac{2^*-q}{2^*}\|u\|_{X_0}^q\nonumber\\
        &\geqslant & \|u\|_{X_0}^2\left( \frac{m_{k+1}^\lambda}{2}-\frac{\varepsilon C}\theta|\Omega|^\frac{2^*-2}{2^*}  -\frac{\delta c}\theta |\Omega|^\frac{2^*-q}{2^*}\|u\|_{X_0}^{q-2} \right)\label{J>0} .
\end{eqnarray}

Taking $\varepsilon >0$ such that $\dfrac{m_{k+1}^\lambda}{4}>\varepsilon\frac{ c}{\theta} |\Omega|^\frac{2^*-2}{2^*}$, we obtain from \eqref{J>0} that 
\begin{equation}
    J_\lambda (u) \geqslant \|u\|_{X_0}^2\left( \frac{m_{k+1}^\lambda}{4}  -\frac{\delta C}\theta |\Omega|^\frac{2^*-q}{2^*}\|u\|_{X_0}^{q-2} \right).
    \label{lastJ32}
\end{equation}
Therefore, for any $u\in X_0^+ $ such that \begin{equation}
 0<   \|u\|_{X_0}\leqslant \rho:=\left(\frac{m_{k+1}^\lambda\theta}{8\delta C |\Omega|^\frac{2^*-q}{2^*}}\right)^\frac{1}{q-2},
\end{equation}
we obtain from (\ref{lastJ32}) that $J_\lambda(u)> 0.$

Thus, we have shown that $J_\lambda$ satisfies the conditions for a local linking at the origin as in (\ref{locallinkcond}).  It then follows from a result due to Liu in \cite[Theorem $2.1$]{LiuS} that the $k_o^{\mbox{th}}$  critical group of $J_\lambda$ at the origin is nontrivial; namely,
\begin{equation}
    C_{k_o}(J_\lambda,0)\not\cong 0,
    \nonumber 
\end{equation}
where $k_o=\dim X_{0}^{-}.$ This concludes the proof of the lemma.
\end{proof}

\section{Critical Groups at Infinity}\label{cinfitysection}

In this section, we will compute the critical groups of $J_\lambda$ at infinity, which were first introduced by Bartsch and Li in \cite{BLi}. 

Let $X$ be a Banach space and assume that $J\in C^1(X,\R)$ satisfies the Palais-Smale condition. Let $\mathcal{K}=\{u\in X\, |\, J^{\prime}(u)=0\}$ be the set of critical points of $J$ and assume that under these assumptions the critical value set is bounded from below; that is, 
$$a_o < \inf J (\mathcal{K}),$$
for some $a_o\in\R.$

The critical groups at infinity are defined by
\begin{equation} 
C_{k}(J ,\infty) = H_{k}(X,J^{a_{o}}),\quad\mbox{ for all }k\in\Z, 
\label{cinfdef}
\end{equation}
(see \cite{BLi}).  These critical groups are well-defined as a consequence of the second deformation theorem (see Perera and Schechter \cite[Lemma $1.3.7$]{PerSch}). 

There are many techniques presented in the literature to address the computation of the critical groups at infinity. One of them is to display a homotopy equivalence between two topological spaces where the homology groups of one the spaces are known. 

Following this line of reasoning, we will show that we can deform a particular sub-level set $J_\lambda^{-M}$, for some $M > 0$, to the unit sphere $S^{\infty}$ in $X_0$. The space $S^{\infty}$ is contractible (see Benyamini {\it et al.} \cite{BS1}) and that will allow us to show that the reduced homology groups of $J^{-M}$ are trivial.  Then, using an argument involving the long exact sequence of the topological pair $(X_0,J_\lambda^{-M})$, we will achieve our goal.

\begin{lemma}
Let $K:\R^{N}\backslash\{0\}\rightarrow (0,\infty)$ satisfy the assumptions $(H_1)$--$(H_3)$.
    Assume that the hypotheses $(H_1)$--$(H_8)$ are satisfied. Then, there exists $\widetilde{M}>0$ such that, for all $M  \geqslant \widetilde{M}$, $J_\lambda^{-M}$ is homotopically equivalent to $S^{\infty}=\{u\in X_0\, |\, \|u\|_{X_0}=1\}$, the unit sphere in $X_0$.  
\label{lemacinf0}
\end{lemma}
\begin{proof}[Proof:]
It follows from the definition of $J_\lambda$ in (\ref{func0}) that   
\begin{eqnarray}\label{J_to_infty2}
   J_\lambda (tu)&=& \frac{t^2}{2}-\lambda \frac{t^2}{2}\|u\|_{L^2(\Omega)}^2-\int_{\Omega}F(x,tu)\,dx, \quad\mbox{ for }u\in S^{\infty};
   \end{eqnarray}
so that, using the estimate in (\ref{estARbigF}), we get
\begin{eqnarray}\label{J_to_infty}
   J_\lambda (tu)&\leqslant& \frac{t^2}{2}-\lambda \frac{t^2}{2}\|u\|_2^2-a_1t^\mu \|u\|_\mu^\mu+a_2|\Omega|,\quad\mbox{ for }u\in S^{\infty}.
   \end{eqnarray}
Then, since $\mu >2$ and $a_1>0$, $\displaystyle\lim_{t\to \infty} J_\lambda (tu)=-\infty$, for $u\in S^{\infty}$. 

Thus, given $\widetilde{M} > 0$, for a given 
$u\in S^{\infty}$,
there exists $t_0 > 0$ such that, 
\begin{equation}
    J_\lambda(tu) < -\widetilde{M},\quad\mbox{ for all }t>t_{0}.
    \label{condinf0}
\end{equation}

Next, for $u\in S^{\infty}$, 
\begin{equation}
    J_{\lambda}(u) = \frac{1}{2}-\frac{\lambda}{2}\|u\|_{L^2(\Omega)}^2-\int_{\Omega}F(x,u)\,dx, 
    \label{eineq55}
\end{equation}
where, using the estimate in (\ref{estbigF}), 
$$F(x,u)\leqslant \varepsilon\int_{\Omega}|u|^2\,dx + \delta(\varepsilon)\int_{\Omega}|u|^q\,dx.$$
Taking $\varepsilon=1$ in the previous inequality, we get 
\begin{equation}
    \int_{\Omega}F(x,u)\,dx \leqslant \|u\|_{L^2(\Omega)}^2 + \delta_1\|u\|_{L^q(\Omega)}^q, 
    \label{eineq66}
\end{equation}
for some positive constant $\delta_1$. 

Then, using the estimates in (\ref{eineq2}), we obtain from (\ref{eineq66}) that 
\begin{equation}
    \int_{\Omega}F(x,u)\,dx \leqslant (C(\theta))^2 + \delta_1(C(\theta))^{q},\quad\mbox{ for }u\in S^{\infty}.
    \label{eineq56}
\end{equation}

Consequently, we get from (\ref{J_to_infty2}) that 
\begin{equation}\label{J>-M}
\begin{aligned}
    J_{\lambda}(u)&\geqslant \frac{1}{2}-\left(\frac{\lambda}{2}(C(\theta)^2+(C(\theta))^2+\delta_1(C(\theta))^{q}\right)\\
    &\geqslant M_1,
\end{aligned}
\end{equation}
where $M_1$ is a constant. 

Next, observe that 
\begin{equation}
   \frac{1}2t^2 \|u\|_{X_0}^2-\lambda \frac{1}2t^2 \|u\|_{2}^2= J_\lambda (tu)+\lambda \int_\Omega F(x, tu) \, dx.
   \label{jj112}
\end{equation}
Thus, we obtain that
\begin{equation}
\begin{aligned}
\frac{d J_\lambda}{dt}(tu)     &  = t \|u\|_{X_0}^2-\lambda t \|u\|_{2}^2-\lambda \int_\Omega f(x, tu)u \, dx  \\
     &= \frac{2}{t} \left( t^2 \|u\|_{X_0}^2-\lambda t^2 \|u\|_{2}^2- \frac{\lambda}{2} \int_\Omega f(x, tu)tu \, dx \right).
\end{aligned}
\label{jj113}
\end{equation}
Next, substitute (\ref{jj112}) into (\ref{jj113}) and apply the hypothesis $(H_8)$ to get 
\begin{equation}
\begin{aligned}
    \frac{d J_\lambda}{dt}(tu)  
     &= \frac{2}{t} \left( J_\lambda (tu)+\lambda \int_\Omega F(x,tu)\, dx -\frac{\lambda}{2}\int_\Omega f(x, tu)tu \, dx \right)\\
     &\leqslant  \frac{2}{t} \left( -M+\frac{\lambda}\mu \int_\Omega f(x,tu)tu\, dx -\frac{\lambda}{2}\int_\Omega f(x, tu)tu \, dx  +\frac{\lambda a |\Omega|}{\mu}\right)\\
     &\leqslant \frac{2}{t} \left( -M-\lambda \int_\Omega \left(\frac{1}{2}-\frac{1}{\mu}\right) f(x, tu)tu \, dx +\frac{\lambda a |\Omega|}{\mu} \right) \\
     &\leqslant \frac{2}{t} \left( -M+\frac{\lambda a |\Omega|}{\mu} \right),
\end{aligned}
\label{dertest1}
\end{equation}
for $tu\in J_\lambda ^{-M}$, and $t>t_0.$

Thus, choosing $M\geqslant\widetilde{M}$ sufficiently large in (\ref{dertest1}), we get 
\begin{equation} 
\frac{d J_\lambda}{dt}(tu) < 0, \quad\mbox {for all } t>t_0\mbox{ and } u\in S^{\infty}. 
\label{dertest}
\end{equation}

Hence, by virtue of (\ref{J>-M}), \eqref{J_to_infty}, and (\ref{dertest}), it follows from the intermediate value theorem that, for any $u\in S^\infty$, there exists $T=T(u)>1$ such that 
$$J_\lambda (Tu)=-M.$$

We can also conclude from the implicit function theorem (see \cite[Theorem $15.1$]{Deim1}) that $T\in C(S^{\infty},\R)$.

Finally, let $B^{\infty}=\{u\in X_0:\|u\|_{X_0}\leqslant 1\}$, the unit ball in $X_0$. Define the homotopy $\eta:[0,1]\times (X_0\backslash B^{\infty})\rightarrow X_0\backslash B^{\infty}$ by
\begin{eqnarray}
    \eta(t,u) = (1-t)u + tT(u)u,
    \nonumber 
\end{eqnarray}
for $t\in[0,1]$ and $u\in X_0\backslash B^{\infty}$. Observe that $\eta(0,u)=u$ and $\eta(1,u)\in J_\lambda^{-M}.$ Thus, $\eta$ is a deformation retract from $X_0\backslash B^{\infty}$ onto $J_{\lambda}^{-M}$. Since $X_0\backslash B^{\infty}$ is homotopically equivalent to $ S^{\infty}$, we conclude that 
$$J_\lambda^{-M}\cong X_0\backslash B^{\infty}\cong S^{\infty};$$
that is, $J_{\lambda}^{-M}$ is homotopically equivalent to $S^{\infty}$. This concludes the proof of the lemma. 
\end{proof}

Using the fact the space $X_0$ is contractible and Lemma \ref{lemacinf0},  we will conclude that the critical groups of $J_\lambda$ at infinity are all trivial. We present the proof of this claim in the next lemma for the reader's convenience. 

\begin{lemma}
    Let $J_\lambda$ be the functional defined in (\ref{func0}) and assume that the hypotheses $(H_1)$--$(H_8)$ are satisfied. Then, the critical groups of $J_\lambda$ at infinity are given by 
    \begin{equation}
        C_{k}(J_\lambda,\infty)\cong 0,\quad\mbox{ for all }k\in\Z. 
        \label{cinftygroups}
    \end{equation}
\end{lemma}
\begin{proof}[Proof:]
    Consider the topological pair $(X_0,J_\lambda^{-M})$ and its long exact sequence of reduced homology groups given by 
\begin{equation}
\ldots{\rightarrow}{\widetilde H}_{k}(J_\lambda^{-M})\stackrel{i_{*}}{\rightarrow}{\widetilde H}_{k}(X_0)\stackrel{j_{*}}\rightarrow {\widetilde H}_{k}(X_0,J_\lambda^{-M})\stackrel{\partial_{*}}\rightarrow {\widetilde H}_{k-1}(J_\lambda^{-M})\stackrel{i_{*}}\rightarrow\ldots,
\label{longseq11}
\end{equation}
where $i_{*}$ and $j_{*}$ are the induced homomorphisms of the inclusion maps
$$i:J_\lambda^{-M}\rightarrow X_0,\quad j:(X_0,\emptyset)\rightarrow (X_0,J_\lambda^{-M}),$$ respectively, and $\partial_{*}$ is a boundary homomorphism as defined in Hatcher \cite[Page $117$]{AH}. 

Since $X_0$ is contractible and $J_{\lambda}^{-M}$ is homotopically equivalent to $S^{\infty}$ per Lemma \ref{lemacinf0}, which is also contractible, the reduced singular homology groups of $X_0$ and $J_\lambda^{-M}$ are given by 
\begin{equation}
    \widetilde{H}_{k}(X_0)\cong 0 \quad\mbox{ and }\quad \widetilde{H}_{k}(J_\lambda^{-M})\cong 0;\quad\mbox{ for all }k\in\Z,
    \label{defhh}
\end{equation}
respectively. 

Next, substitute (\ref{defhh}) into (\ref{longseq11}) to get  
\begin{equation}
    \ldots{\rightarrow}0\stackrel{i_{*}}{\rightarrow}0\stackrel{j_{*}}\rightarrow {\widetilde H}_{k}(X_0,J_\lambda^{-M})\stackrel{\partial_{*}}\rightarrow 0 \stackrel{i_{*}}\rightarrow\ldots.
\label{longseq01}
\end{equation}

Hence, it follows by the exactness of the sequence (\ref{longseq01}) that 
 $$C_{k}(J_\lambda,\infty)=\widetilde{H}_{k}(X_0,J_\lambda^{-M})\cong 0,\quad\mbox{ for all }k\in\Z.$$

\end{proof}

\section{Existence of Multiple Solutions}\label{secsolsection}

In this section, we prove the existence of at least three nontrivial solutions of problem (\ref{prob_1}) for the case $\lambda < \lambda_1$ and the existence of at least two nontrivial solutions for the case $\lambda\geqslant\lambda_1$, respectively.

In \cite[Corollary $21$]{Ser3}, Servadei and Valdinoci used a cutoff technique to prove that problem (\ref{prob_1}) admits a non-negative solution $u_{+}\in X_{0}$ and a non-positive solution $u_{-} \in X_0$; both solutions are given by the mountain pass theorem of Ambrosetti-Rabinowitz (see \cite[Corollary $21$]{Ser3}) for the case $\lambda<\lambda_1$. We will present an extension of this result in the next theorem using a line of reasoning similar to that in \cite[Section $3$]{Wang}. 

\begin{theorem}
    Let all the assumptions of Theorem \ref{sertheo1} be satisfied. Then, for any $\lambda < \lambda_1$, problem (\ref{prob_1}) admits at least three nontrivial solutions where one is non-positive and another is non-negative. 
    \label{maintheo2}
\end{theorem}

\begin{proof}[Proof:]
Assume that the critical set, $\mathcal{K}$, of $J_\lambda$ has only three critical points; namely, $\mathcal{K}=\{0,u_{+},u_{-}\}$, where $u_{+}$ and $u_{-}$ are the critical points obtained via the   mountain-pass theorem in Servadei and Valdinoci \cite[Corollary 21]{Ser3}.  Let $c_1=J_\lambda(u_+)$ and $c_{2}=J_\lambda(u_{-})$ be the critical values of $J_\lambda$ at $u_+$ and at $u_{-}$, respectively. 
It also follows from the result in 
\cite[Corollary 21]{Ser3} 
that $c_1>0$ and $c_2>0$.  We may assume, without loss of generality, that $0<c_2\leqslant c_1$.

Let $a_1,a_2,a_3\in\R$ be such that
    \begin{equation}
        a_1 < 0 < a_2 < c_2\leqslant  c_1 < a_3.
        \label{cvcond1}
    \end{equation}

By virtue of Proposition \ref{propMMP}, we have
\begin{equation}
    C_{1}(J_\lambda,u_+)\not\cong 0\quad\mbox{ or}\quad C_{1}(J_\lambda,u_{-})\not\cong 0,
    \label{cgroupsposneg}
\end{equation}
since $u_{\pm}$ are the critical points given by the mountain-pass theorem found in \cite[Corollary $21$]{Ser3}.

Consider the triple $(J_\lambda^{a_{1}},J_\lambda^{a_{2}},J_\lambda^{a_{3}})$ and its long exact sequence of reduced singular homology groups given by 
\begin{equation}
\ldots{\rightarrow}\widetilde{H}_{k+1}(J_\lambda^{a_{2}},J_\lambda^{a_{1}})\stackrel{i_{*}}{\rightarrow}\widetilde{H}_{k+1}(J_\lambda^{a_{3}},J_\lambda^{a_{1}})\stackrel{j_{*}}\rightarrow\widetilde{H}_{k+1}(J_\lambda^{a_{3}},J_\lambda^{a_{2}})\stackrel{\partial_{*}}\rightarrow \widetilde{H}_{k}(J_\lambda^{a_{2}},J_\lambda^{a_{1}})\stackrel{i_{*}}\rightarrow,
\label{longseq211}
\end{equation}
(see \cite[Section $2.1$, page 118]{AH}).

By virtue of Proposition \ref{mainprop}, we have 
\begin{equation}
    \widetilde{H}_{k+1}(J_\lambda^{a_3},J_\lambda^{a_2})\cong C_{k+1}(J_\lambda,u_+)\oplus C_{k+1}(J_\lambda,u_{-}),\quad\mbox{ for all }k\in\Z.
    \label{iso211}
\end{equation}

Substitute (\ref{iso211}) and (\ref{cinftygroups}) into (\ref{longseq211}) to get 
 \begin{equation}
\ldots{\rightarrow}0\stackrel{j_{*}}\rightarrow  C_{k+1}(J_\lambda,u_+)\oplus C_{k+1}(J,u_{-})\stackrel{\partial_{*}}\rightarrow \widetilde{H}_{k}(J_\lambda^{a_{2}},J_\lambda^{a_{1}})\stackrel{i_{*}}\rightarrow 0\stackrel{j_{*}}\rightarrow\ldots,
\label{longseq311}
\end{equation}
for all $k\geqslant 0.$

Using the fact that the sequence (\ref{longseq311}) is exact, we can show that $\partial_*$ is an isomorphism; that is, 
\begin{equation} \widetilde{H}_{k}(J_\lambda^{a_{2}},J_\lambda^{a_{1}})\cong C_{k+1}(J_\lambda,u_+)\oplus C_{k+1}(J,u_{-}),\quad\mbox{ for all }k\in \Z.
\label{iso99}
\end{equation}

By virtue of Proposition \ref{mainprop}, Lemma  (\ref{lemmacorigincase1}) and (\ref{cvcond1}), we have 
\begin{equation}
     H_{k}(J_\lambda^{a_2},J_{\lambda}^{a_1}) \cong C_{k}(J_\lambda,0) \cong \delta_{k,0}\F,\quad\mbox{ for } k\in\Z. 
    \label{interc1}
\end{equation}
In particular, setting $k=0$ in (\ref{interc1}) and using the definition of the reduced homology groups as in (\ref{redhomrel}), we conclude that 
\begin{equation}
    \widetilde{H}_{0}(J_\lambda^{a_2},J_\lambda^{a_1})\cong 0.
    \label{interc2}
\end{equation}
Thus, setting $k=0$ in (\ref{iso99}) and using the assertions (\ref{cgroupsposneg}) and (\ref{interc2}), we get  
\begin{equation}
    0\cong \widetilde{H}_{0}(J_\lambda^{a_{2}},J_\lambda^{a_{1}})\cong C_{1}(J_\lambda,u_+)\oplus C_{1}(J_\lambda,u_{-})\not\cong 0,
    \nonumber 
\end{equation}
which is a contradiction.  

Hence, the critical set $\mathcal{K}$ must have at least three nontrivial critical points. This concludes the proof of the theorem. 
\end{proof}

The next theorem presents the existence of a second nontrivial solution of problem (\ref{prob_1}) for the case $\lambda\geqslant \lambda_1$. Theorem \ref{cingotheo} will play a key role in the proof of this result. 

\begin{theorem}\label{maintheopaper}
     Let $s\in (0,1), N>2s$ and $\Omega$ be an open bounded set of $\R^{N}$ with Lipschitz boundary. Let $K:\R^{N}\backslash\{0\}\rightarrow (0,+\infty)$ be a function satisfying conditions 
     $(H_1)$--$(H_3)$.  Assume that $f$ satisfies  hypotheses $(H_4)$, $(H_6)$--$(H_9)$.  Let $u_o$ be the solution found in Theorem \ref{sertheo1} through the linking theorem. Assume also that 
     $\lambda\in [\lambda_{k},\lambda_{k+1})$, for some  $k\in\N$, and define $X_0^{-}$ as in 
     (\ref{hilbdecomp}).  Then,      
     problem (\ref{prob_1}) admits at least two nontrivial solutions provided that
         \begin{equation}
         k_o+1\not\in[\mu_o,\mu_o+\nu_o],
         \label{kocond}
         \end{equation}
         where $k_o=\dim X_{0}^{-}$, $\mu_o$ and $\nu_0$ are the Morse index and nullity of $u_o$, respectively. 
  
\end{theorem}
\begin{proof}[Proof:]
Let $\lambda\in [\lambda_k,\lambda_{k+1})$, for some $k\in\N$, and $u_o$ be the critical point of $J_\lambda$ found in \cite[Section $4.2$]{Ser3}. The authors of \cite{Ser3} proved that the functional $J_\lambda$ has the geometry of the the linking theorem of Rabinowitz \cite{Rab78}; they also proved that any Palais-Smale sequence in $X_0$ is bounded (see  \cite[Proposition $20$]{Ser3}) and possesses a convergent subsequence (see \cite[Proposition $14$]{Ser3}). Hence $J_\lambda$ satisfies the Palais-Smale condition.

Assume, by a way of contradiction, that the critical set of $J_\lambda$ contains only two critical points; namely, $\mathcal{K}=\{0,u_o\}$.  Let $c=J_\lambda(u_o)$ be the critical value of $J_\lambda$ at $u_o$.  It follows from the results
in \cite{Ser3} that $c>0$ (see 
\cite[Proposition 17, page 2120]{Ser3}).
Let $a_1,a_2,a_3\in\R$ be such that
    \begin{equation}
        a_1 < 0 < a_2 < c < a_3.
        \label{cvcond}
    \end{equation}
Consider the triple $(J_\lambda^{a_{1}},J_\lambda^{a_{2}},J_\lambda^{a_{3}})$ and its long exact sequence of reduced homology groups given by
\begin{equation}
{\rightarrow}\widetilde{H}_{k+1}(J_\lambda^{a_{2}},J_\lambda^{a_{1}})\stackrel{i_{*}}{\rightarrow}\widetilde{H}_{k+1}(J_\lambda^{a_{3}},J_\lambda^{a_{1}})\stackrel{j_{*}}\rightarrow\widetilde{H}_{k+1}(J_\lambda^{a_{3}},J_\lambda^{a_{2}})\stackrel{\partial_{*}}\rightarrow \widetilde{H}_{k}(J_\lambda^{a_{2}},J_\lambda^{a_{1}})\stackrel{i_{*}}\rightarrow
\label{longseq2}
\end{equation}

Since $J_\lambda$ satisfies the Palais-Smale condition, it follows from (\ref{cinftygroups}) that 
\begin{equation}
    \widetilde{H}_{k+1}(J_\lambda^{a_3},J_\lambda^{a_1})\cong \widetilde{H}_{k+1}(X_0,J_\lambda^{a_1})\cong C_{k+1}(J_\lambda,\infty)\cong 0 ,\quad\mbox{ for }k \in\Z.
    \label{iso1}
\end{equation}

By virtue of Proposition \ref{mainprop}, we have 
\begin{equation}
    \widetilde{H}_{k+1}(J_\lambda^{a_3},J_\lambda^{a_2})\cong C_{k+1}(J_\lambda,u_o),\quad\mbox{ for all }k\in\Z.
    \label{iso2}
\end{equation}

Substitute (\ref{iso1}) and (\ref{iso2}) into (\ref{longseq2}) to get 
 \begin{equation}
\ldots{\rightarrow} 0\stackrel{j_{*}}\rightarrow  C_{k+1}(J_\lambda,u_o)\stackrel{\partial_{*}}\rightarrow \widetilde{H}_{k}(J_\lambda^{a_{2}},J_\lambda^{a_{1}})\stackrel{i_{*}}\rightarrow 0 \stackrel{j_{*}}\rightarrow\ldots,
\label{longseq3}
\end{equation}
for all $k\geqslant 0.$ 

It follows from the exactness of  the sequence (\ref{longseq3}) that $\partial_{*}$ is an isomorphism. In fact, $\ker(\partial_{*})=Im(j_{*})=\{0\}$ implies that $\partial_{*}$ is one-to-one and $Im(\partial_*)=\ker (i_{*})$ implies that $\partial_*$ is onto because $Im(i_*)=\{0\}$. We can then write 
\begin{equation} \widetilde{H}_{k}(J_\lambda^{a_{2}},J_\lambda^{a_{1}})\cong C_{k+1}(J_\lambda,u_o),\quad\mbox{ for all }k\in\Z.\label{iso51}
\end{equation}

By virtue of (\ref{iso51}) and (\ref{cgrouporigincase2}) we get that 
\begin{equation}
C_{k}(J_\lambda,0)\cong \widetilde{H}_{k}(J_\lambda^{a_{2}},J_\lambda^{a_{1}})\cong C_{k+1}(J_\lambda,u_o),\quad\mbox{ for all }k\in\Z.\label{iso53}
\end{equation}

In what follows, we analyze the case in which $u_o$ is a nondegenerate isolated critical point of $J_\lambda$ and the 
case in which $u_o$ is degenerate separately.   

In the case in which $u_o$ is a nondegenerate isolated critical point of $J_\lambda$ with finite Morse index $\mu_o$,  the critical groups of $J_\lambda$ at $u_o$ depend only on the Morse index of $u_o$; namely,  
\begin{equation}
    C_{k}(J_\lambda,u_o)\cong\delta_{k,\mu_o}\F,\quad\mbox{ for all }k\in\Z,
    \label{cnondeg1}
\end{equation}
(see \cite[Corollary $8.3$]{MW1}). 

Substitute (\ref{cnondeg1}) in (\ref{iso53}) to get 
\begin{equation} 
C_{k}(J_\lambda,0)\cong C_{k+1}(J_\lambda,u_o)\cong \delta_{k+1,\mu_o}\F,\quad\mbox{ for all }k\in\Z.
\label{ssd1}
\end{equation}

Hence, by virtue of (\ref{cgrouporigincase2}) and (\ref{ssd1}) with $k=k_o$ we get 
\begin{equation} 
0\not\cong C_{k_o}(J_\lambda,0)\cong \delta_{k_o+1,\mu_o}\F\cong 0,\quad\mbox{ for all }k\in\Z,
\label{ssd2}
\end{equation}
which is a contradiction because $k_o+1\ne\mu_o$ by assumption (\ref{kocond}). 

 Next, assume that $u_o$ is a degenerate critical point for $J_\lambda$ with finite Morse index $\mu_0$ and nullity $\nu_0$. Then, it follows from the shifting theorem \cite[Theorem $8.4$] {MW1}  that 
\begin{equation}
    C_{k}(J_\lambda,u_o)\cong C_{k-\mu_o}(\widetilde{J}_\lambda,u_o),\quad\mbox{ for }k\in\Z,
    \label{shift1}
\end{equation}
 where $\widetilde{J}_\lambda : B_{\rho}(u_o)\rightarrow \R$ is the $C^2$ function in Theorem \ref{cingotheo},  $B_{\rho}(u_o)$ is a small ball in $V:=\ker J_\lambda^{\prime\prime}(u_o)$ centered at $u_o$, and radius $\rho>0$, and $\nu_o=\dim V.$

We will analyze the three possible cases for the critical point $u_o$ as in Theorem \ref{cingotheo} and show that we reach a contradiction in each one of them. 

First, by virtue of the Remark \ref{fredremark}), the operator $J^{\prime\prime}(u_o)$ is a Fredholm operator. Next, if $u_o$ is a local minimum for $\widetilde{J}_\lambda$, then by virtue of (\ref{shift1}) with $k=k_o$, (\ref{iso53}), 
and \cite[Example $1$, page $33$]{KC}, we have 
\begin{equation}
  0\not\cong C_{k_o}(J_\lambda,0)\cong C_{k_o+1}(J_\lambda,u_o)\cong C_{k_o+1-\mu_o}(\widetilde{J}_\lambda,u_o)\cong \delta_{k_o+1-\mu_o,0}\F. 
 \label{lmin}
\end{equation}
 On the other hand, since $k_o+1\ne \mu_o$ by assumption (\ref{kocond}),   
$\delta_{k_o+1-\mu_o,0}\F \cong 0$. This yields a contradiction in view of (\ref{lmin}).

For the case in which $u_o$ is a local maximum for $\widetilde{J}_\lambda$, it follows from (\ref{shift1}), (\ref{iso53}), and \cite[Example $1$, page $33$]{KC} that 
\begin{equation*}
  0\not\cong C_{k_o}(J_\lambda,0)\cong C_{k_o+1}(J_\lambda,u_o)\cong C_{k_o+1-\mu_o}(\widetilde{J}_\lambda,u_o)\cong \delta_{k_o+1-\mu_o,\nu_o}\F , 
\end{equation*}
which leads to a contradiction since $k_o+1\ne \mu_o+\nu_o$ by the assumption (\ref{kocond}).

It remains to consider the case in which $u_o$ is neither a local maximum or a local minimum for $\widetilde{J}_\lambda$.  In this case, combine 
(\ref{iso53}) and (\ref{cgrouporigincase2}) to 
obtain
 \begin{equation} 
0\not\cong C_{k_o}(J_\lambda,0)\cong C_{k_o+1}(J_\lambda,u_o).
 \label{nminmax}
 \end{equation}
Since we are assuming that $u_o$ is neither a local maximum or a local minimum for $\widetilde{J}_\lambda$, we can use part (c) of
Theorem \ref{cingotheo}, and assumption (\ref{kocond}), to get
\begin{equation}C_{ko+1}(J_\lambda,u_o)\cong 0.
\label{finaleq}
\end{equation} 
Finally, substituting (\ref{finaleq}) into (\ref{nminmax}) leads to a  contradiction. This concludes the proof of the theorem. 
 \end{proof}

\section{Second Multiplicity Result}\label{secexistsec}

In this section, we will prove that the functional $J_\lambda$, for $\lambda$ between two consecutive
eigenvalues, $\lambda_k$ and $\lambda_{k+1}$,
has an unbounded sequence of critical points when the non-linearity $f(x,s)$ is odd in the second variable $s$. The result will follow as an application of the $\Z_2$ version of the Mountain Pass Theorem of Ambrosetti-Rabinowitz \cite{AmbRab}. For the reader's convenience, we present the version of this theorem given in \cite[Theorem $9.12$]{Rabinowitz1}. 

\begin{theorem}\cite[Theorem $9.12$]{Rabinowitz1}
Let $E$ be an infinite-dimensional Banach space, and let $J\in C^{1}(E,\R)$ be even. Suppose that $J$ satisfies the Palais-Smale condition, and $J(0)=0.$ Assume also that $E=V\oplus W,$ where $V$ is finite dimensional, and $J$ satisfies
\begin{itemize}
    \item [(i)] there are constants $\rho,\alpha>0$ such that $J|_{\partial B_\rho\cap W}\geqslant\alpha,$ and 
    \item [(ii)] for each finite dimensional subspace $\widetilde{E}\subset E$, there is an $R=R(\widetilde{E})>0$ such that $J\leqslant 0$ on $\widetilde{E}\backslash B_{R(\widetilde{E})}$.
\end{itemize}
Then $J$ has an unbounded sequence of critical values. 
\label{z2theo}
\end{theorem}

To establish the second multiplicity result for problem (\ref{prob_1}), we will consider the decomposition of the Hilbert space $X_0=V\oplus W$ with $V=X_{0}^{-}$ and $W=X_{0}^{+}$, where $X_0^{-}$ is the finite-dimensional space defined in (\ref{hilbdecomp}). 

\begin{theorem}
     Suppose that 
  $K:\R^{N}\backslash\{0\}\rightarrow (0,+\infty)$  satisfies conditions $(H_1)$--$(H_3)$, and $f$ satisfies $(H_4),$ $(H_5)$, $(H_8)$ and $f(x,s)$ is odd in $s$. Assume also that $\lambda\in [\lambda_{k},\lambda_{k+1})$, for some $k\in\N$, and $X=V\oplus W$ with $V=X_{0}^{-}$ and $W=V^{\perp}$ as in (\ref{hilbdecomp}). Then, problem (\ref{prob_1}) possesses an unbounded sequence of weak solutions.
    \label{secexistheo}
\end{theorem}
\begin{proof}[Proof:] We will show that the functional $J_\lambda$, defined in (\ref{func0}), satisfies the conditions of Theorem \ref{z2theo}. 

First, recall that 
\begin{equation}
    J_\lambda(u)=\frac{1}{2}\|u\|_{X_0}^2 -\frac{\lambda}{2} \int_{\Omega} u^2 dx-\int_\Omega F(x,u)dx, \quad\mbox{ for }u\in X_{0}.
    \label{func01}
\end{equation}
It was shown in \cite{Ser3} that $J_\lambda\in C^{1}(X_0,\R)$ and $J_\lambda$ also satisfies the Palais-Smale condition for  $\lambda\geqslant \lambda_1$ in \cite[Propositions $14$,$20$]{Ser3}.  By virtue of the fact that $f(x,s)$ is odd in $s$, we have that $J_\lambda$ is even.  We also have that $J_\lambda(0)=0$. Hence, it remains to verify that conditions $(i)$ and $(ii)$ in Theorem \ref{z2theo} hold true. 

First, we prove that condition $(i)$ is satisfied. 

For every  $\varepsilon >0$, there exists a $\delta(\varepsilon)>0$ such that  
 \begin{equation}
    \int_{\Omega}F(x,u)\,dx   
       \leqslant \left(C_{1}\varepsilon + C_2\delta(\varepsilon)\|u\|_{X_{0}}^{q-2}\right)\|u\|_{X_0}^2,
       \label{ineq666}
\end{equation}
for $u\in X_0$, where $q\in (2,2^{*})$, $\theta > 0$, $C_1,C_2$ are positive constants as was proved in (\ref{estbigF}), which is derived from (\ref{estlittlef}). 

Let $u\in W$. It follows from the characterization of the eigenvalue $\lambda_{k+1}$ in \cite[Proposition $9$]{Ser3} that 
\begin{equation}
    \lambda_{k+1}\int_{\Omega}|u|^2\,dx\leqslant \|u\|_{X_0}^2,
    \quad\mbox{ for all }u\in W.
    \label{chalk1}
\end{equation}

Then, by virtue of (\ref{func01}), (\ref{ineq666}), and (\ref{chalk1}), we obtain from (\ref{func01}) that 
\begin{equation}
    J_\lambda(u) \geqslant \left[\frac{1}{2}-\frac{\lambda}{2\lambda_{k+1}}- C_{1}\varepsilon  - C_2\delta(\varepsilon)\|u\|_{X_{0}}^{q-2}\right]\|u\|_{X_0}^2,\quad\mbox{ for }u\in W.
    \label{ineq1k}
\end{equation}

Let $\|u\|_{X_0}=\rho.$   Then, we can rewrite (\ref{ineq1k}) as 
\begin{equation}
    J_\lambda(u) \geqslant \left[\frac{1}{2}-\frac{\lambda}{2\lambda_{k+1}}- C_{1}\varepsilon   - C_2\delta(\varepsilon)\rho^{q-2}\right]\rho^2,\quad\mbox{ for }u\in W.
    \label{ineq2k}
\end{equation}

We next use the assumption that $\lambda_{k}\leqslant\lambda < \lambda_{k+1}$, for some $k\in\N$. 

Choose $\rho$ in (\ref{ineq2k}) such that 
\begin{equation}
    \frac{1}{2}-\frac{\lambda}{2\lambda_{k+1}}-C_{1}\varepsilon  - C_2\delta(\varepsilon)\rho^{q-2}\geqslant\frac{1}{4}.
    \label{ine888}
\end{equation}
This choice of $\rho>0$ is possible since we can choose $\varepsilon=\dfrac{\theta \lambda}{4C_1\lambda_{k+1}} $ in the estimate (\ref{ineq666}). Then, it follows from (\ref{ineq2k}) and the choice of $\rho$ in (\ref{ine888}) that 
\begin{equation}
    J_\lambda(u)\geqslant \frac{1}{4}\rho^2 :=\alpha,\quad\mbox{ for }u\in\partial B_{p}\cap W,
    \nonumber 
\end{equation}
where $\partial B_{\rho}=\{u\in X_0 :  \|u\|_{X_0}=\rho\}$. This establishes condition $(i)$ in Theorem \ref{z2theo}. 

Next, we proceed to check condition $(ii)$ in Theorem (\ref{z2theo}).
 
Let $\widetilde{E}\subset X_0$ be a finite-dimensional space such that $\widetilde{E}=span\{e_1, \dots , e_n\},$ where $\{e_i\}_{i=1}^{n}$ is an orthonormal subset of $X_0$.  
Set 
$$\|u\|_{\widetilde{E}}:=\max_{i\in [1,n]}|u_i|,$$
where $u=u_1e_1+\dots+u_ne_n$, $u_i\in \mathbb{R}$. 
We note that $\|\cdot\|_{\widetilde{E}}$ defines a norm in $\widetilde{E}$, which is equivalent to $\|\cdot \|_{X_0}$ and $\|\cdot \|_{L^\mu(\Omega)}$, 
when restricted to the finite dimensional subspace
$\widetilde{E}$. Thus, there exist $\hat C_1, \hat C_2 >0$ such that, for all $u\in \widetilde{E}$, we get 
\begin{equation}\label{equi_norms}
 \|u\|_{X_0} \leq \hat C_1 \|u\|_{\widetilde{E}} \quad \text{and} \quad \hat C_2   \|u\|_{\widetilde{E}} \leqslant \|u\|_{L^\mu(\Omega)}.   
\end{equation}
Since $\mu >2$, there exists $R > 0$ such that, for
$t\geqslant R$, we obtain
\begin{equation}
 \frac{t^2}{2}\hat C^2_1  -a_1 \hat C^\mu _2 t^\mu +a_2|\Omega| \leqslant 0. \label{Je_i<0} 
\end{equation}
Let $u\in \widetilde{E}$ with $\|u\|_{\widetilde{E}}\geqslant R$. Write $u=u_1e_1+\cdots + u_ne_n$, $u_i\in \mathbb{R}$ and
choose $i\in \{1,2,3,\ldots,n\}$ so that 
$|u_i|=\|u\|_{\widetilde{E}}$.   Use the definition of $J_\lambda$ in  (\ref{func01}) and the
superquadratic growth condition in (\ref{estARbigF}) 
to obtain the estimate
$$
 J_\lambda (u) \leqslant \frac{1}{2}\|u\|_{X_0}^2-\|u\|^\mu_{L^\mu(\Omega)} +a_2|\Omega|,
    \quad\mbox{ for } u\in X_0;
$$
so that, using the inequalities in \eqref{equi_norms},
$$
 J_\lambda (u) \leqslant \frac{1}{2} \hat C^2_1\|u\|_{\widetilde{E}}^2-a_1\hat C_2^\mu\|u\|^\mu_{\widetilde{E}}+a_2|\Omega|,
    \quad\mbox{ for }u\in \widetilde{E}.
$$
It then follows from \eqref{Je_i<0}  that 
$$
    \begin{aligned}
        J_\lambda (u) &\leqslant  
        \frac{1}{2} \hat C^2_1|u_i|^2-a_1\hat C_2^\mu|u_i|^\mu+a_2|\Omega|        \leqslant 0,
            \quad\mbox{ for } u\in \widetilde{E}
                \mbox{ with } \| u\|_{\widetilde{E}}
                    \geqslant R.
    \end{aligned}
$$
This establishes condition $(ii)$ of Theorem \ref{z2theo} by virtue of the equivalence of the norms in $\widetilde{E}$.

Therefore, since all the hypotheses of Theorem \ref{z2theo} are satisfied, problem (\ref{prob_1}) has a sequence of unbounded critical values $c_k=J_{\lambda}(u_k)$, where $u_k$ is a weak solution of problem (\ref{prob_1}), for $k\in\N.$

Next, we show that the sequence $(u_k)_{k\in\N}$ is unbounded in $X_0$. In fact, since $u_k$ is a weak solution for problem (\ref{prob_1}), it follows from (\ref{frechet0}) with $\varphi=u_k$ that 
\begin{equation}
    \langle J_\lambda^\prime(u_k),u_k\rangle =0; \quad\mbox{ for all }k\in\N,
    \nonumber
\end{equation}
that is, 
\begin{equation}
    \|u_k\|_{X_0}^2 = \lambda\int_{\Omega}u_k^2\,dx + \int_{\Omega}f(x,u_k)u_k\,dx,
    \label{ffrr11}
\end{equation}
for all $k\in\N.$

Substitute (\ref{ffrr11}) into (\ref{func01}) to get
\begin{equation}
    c_{k}=J_\lambda(u_k)=\frac{1}{2}\int_{\Omega}\left(f(x,u_k)u_k-2F(x,u_k)\right)\,dx,
        \quad \mbox{ for all } k\in\N.
    \label{ffrr12}
\end{equation}

Use (\ref{ffrr11}) to obtain the identity 
$$
     \|u_k\|_{X_0}^2 = \int_{\Omega}\left(f(x,u_k)u_k-2F(x,u_k)\right)\,dx + \int_{\Omega}2F(x,u_k)\,dx.
        \quad\mbox{ for } k\in\N,
$$
from which we get the inequality
\begin{equation}\label{AReqo05}
    \|u_k\|_{X_0}^2 \geqslant \int_{\Omega}\left(f(x,u_k)u_k-2F(x,u_k)\right)\,dx + 2\|u_k\|_{L^{\mu}(\Omega)}^{\mu}-C,
\end{equation}
for $k\in\N$ and $C$ is a positive constant,
where we have used the estimate in (\ref{estARbigF}).

It follows from (\ref{AReqo05}) and (\ref{ffrr12}) 
that
\begin{equation}
    \begin{aligned}
        \|u_k\|_{X_0}^2  & \geqslant 
            2c_k - C,
                \quad\mbox{ for } k\in\N.
    \end{aligned}
    \label{final1}
\end{equation}
The estimate in (\ref{final1}), in conjunction with
the fact that the sequence $(c_k)_{k\in\N}$ is 
unbounded, implies that the sequence $(u_k)_{k\in\N}$ of weak solutions of problem (\ref{prob_1}) is unbounded in $X_0$. 
\end{proof}

\bibliographystyle{plain}
\bibliography{main}

\end{document}